\documentstyle[twoside,leqno]{aop}
\newcommand{\bea}{\begin{eqnarray}}
\newcommand{\ba}{\begin{array}}
\newcommand{\bean}{\begin{eqnarray*}}
\newcommand{\ea}{\end{array}}
\newcommand{\eea}{\end{eqnarray}}
\newcommand{\eean}{\end{eqnarray*}}
\newcommand{\be}{\begin{equation}}
\newcommand{\ee}{\end{equation}}

\newcommand{\lab}{\label}

\def\qed{$\hfill \Box$}
\def\deq{\, {\stackrel {def} {=}}}

\def\sif{\sigma \mbox{\rm -field}}

\def\ep{\epsilon}
\def\E{{\bf{E}}}
\def\P{{\bf{P}}}
\def\N{{\bf N}}

\def\Z{{\bf{Z}}}

\def\G{{\mathcal{G}}}

\def\mc{\Phi}

\def\rF{{\cal S}}

\def\ren{\Upsilon}

\def\measn{\alpha}
\def\measa{\beta}
\def\measb{\pi}

\def\R{{\bf R}}

\def\param{\varphi}
\def\whpsi{\widehat{\psi}}

\def\|{\, |\! | \, }
\def\one{{\bf 1}}

\def\phi{\varphi}
\def\tail{{\cal T}}
\def\and{\, \mbox{ and } \,}
\def\given{\; | \; }
\def\iid{i.i.d.\ }
\def\deq{\, {\stackrel {{\rm def}} {=}} \,}
\def\b2b{\overline{\beta^2}}
\def\var{{\tt Var}}
\def\cov{{\tt Cov}}

\def\shift{\Theta}

\def\coins{\{-1,1\}}
\def\bin{\{0,1\}}
\def\bias{\theta}
\def\renp{\Delta}
\def\rent{\delta}
\def\rhs{right--hand side }
\def\np{Y}
\def\nt{y}
\def\legt{\Lambda^{*}}
\def\sp{Z}
\def\st{z}

\def\ma{\mu_{\measa}}
\def\HH{{\mathcal H}}
\def\mb{\mu_{\measb}}

\def\Q{{\bf Q}}

\def\measu{\eta}

\def\muu{\mu_\eta}
\def\ls{\widehat{R}}
\def\c{\zeta}
\def\h{h}
\def\H{H}
\newtheorem {thm}{Theorem}[section]
\newtheorem {prop}[thm]{Proposition}
\newtheorem {lem}[thm]{Lemma}

\begin{document}


\SPECFNSYMBOL{}{1}{2}{}{}{}{}{}{}%

\AOPMAKETITLE


\AOPAMS{Primary 60G30; secondary 60K35}
\AOPKeywords{Mutually singular measures, Kakutani's dichotomy, renewal 
sequences, random walks, scenery with noise, phase transitions}
\AOPtitle{A Phase Transition in Random Coin Tossing}
\AOPauthor{David A. Levin, Robin Pemantle\footnote{Research
    partially supported by a Presidential Faculty Fellowship} and Yuval
Peres\footnote{Research partially supported by NSF grant \# DMS-9404391}
}
\AOPaffil{University of California,
 University of Wisconsin, and Hebrew University}
\AOPlrh{D.A. Levin, R. Pemantle, and Y. Peres}
\AOPrrh{Random Coin Tossing}
\AOPAbstract{
Suppose
that a coin with bias $\theta$ is tossed at renewal times of a
renewal process, and a fair coin is tossed at all other times.  Let
$\mu_\theta$ be the distribution of the observed sequence of coin tosses,
and let $u_n$ denote the chance of a renewal at time $n$.
Harris and Keane in \cite{HK} showed that if $\sum_{n=1}^\infty u_n^2 = \infty$, then
$\mu_\theta$ and $\mu_0$ are singular,
 while if $\sum_{n=1}^{\infty} u_n^2 < \infty$ and
$\theta$ is small enough, then $\mu_\theta$ is absolutely continuous
with respect to $\mu_0$.  They conjectured
that absolute continuity should not depend on $\theta$, but only on
the square-summability of $\{u_n\}$.
We show that in fact the power law governing the decay
of   $\{u_n\}$ is crucial, and for some renewal sequences $\{u_n\}$,
there is a {\em phase transition\/} at a critical parameter
$\theta_c \in (0,1)$: for
$ |\theta|<\theta_c$ the measures
$\mu_\theta$ and $\mu_0$ are mutually absolutely continuous, but for
$|\theta|>\theta_c$, they are singular.  
We also prove that when $u_n=O(n^{-1})$, the measures
$\mu_\theta$ for $\theta \in [-1,1]$
are all  mutually absolutely continuous.
}

\maketitle

\BACKTONORMALFOOTNOTE{3}

\section{Introduction.} \label{sec:intro}

A {\em coin toss} with {\em bias} $\theta$ is a $\{-1,1\}$-valued
random variable with mean $\theta$, and a {\em fair coin} is a coin
toss with mean zero.  Kakutani's dichotomy for independent sequences
reduces, in the case of coin tosses, to the following:

{\sc Theorem A} (\cite{Ka})
{\em Let $\mu_0$ be the distribution
of \iid fair coin tosses on $\{-1,1\}^{\N}$, and let $\nu_\theta$ be the
distribution of independent coin tosses with biases
$\{\theta_n\}_{n=1}^{\infty}$.
\begin{description}
\item{(i)}
If $\sum_{n=1}^{\infty}\theta_n^2 = \infty$ then $\nu_\theta \perp
\mu_0$, where $\nu \perp \mu$ means that the measures $\nu$ and $\mu$ are
mutually singular.
\item{(ii)}
If $\sum_{n=1}^{\infty}\theta_n^2 <
\infty$, then $\nu_\theta \ll \mu_0$ and $\mu_0 \ll \nu_\bias$,
where $\nu \ll \mu$ means that $\nu$ is absolutely continuous with respect
to $\mu$.
\end{description}}

\medskip
\noindent
For a proof of Theorem A see, for example, Theorem 4.3.5 of \cite{Du}.

Harris and Keane \cite{HK} extended Theorem A(i) to sequences
with a specific type of dependence.  Let $\{\Gamma_n\}$ be a
(hidden) recurrent Markov chain with initial state $o$, called
the {\em origin}.
Suppose that whenever $\Gamma_n=o$,
 an independent coin with bias $\theta \ge 0$ is tossed,
while at all other times an independent fair coin is tossed. Write
$X=(X_1,X_2,\ldots)$ for the record of coin tosses, and let
$\mu_\bias$ be the distribution of $X$.
Let  $\renp_n = \one_{\{\Gamma_n = o\}}$ and denote by
$$
u_n = \P[\Gamma_n = o] =\P[\renp_n=1]
$$
the probability of a return of the chain to the origin at time $n$.
The random variables $\{\renp_n\}$ form a
{\em renewal process}, and their joint distribution is determined
by the corresponding {\em renewal sequence} $\{u_n\}$; see the next section.
Harris and Keane established the following theorem.

{\sc Theorem B} (\cite{HK}) {\em \begin{description}
        \item{(i)} If $\, \sum_{n=1}^\infty u_{n}^{2} = \infty$, then
        $\mu_{\theta} \perp \mu_{0}$.
    \item{(ii)} If  $\, \sum_{n=1}^\infty u_{n}^{2} = \|u\|^2 < \infty$ and
$\theta <  \|u\|^{-1}$,  
        then $\mu_{\bias} \ll\mu_{0}$.
 \end{description}}

Harris and Keane conjectured that singularity of the two laws
$\mu_{\theta}$ and $\mu_{0}$ should not depend on $\theta$, but only
on the return probabilities $\{u_{n}\}$.  In
particular, they asked whether the condition $\sum_{k=0}^\infty u_k^2
< \infty$ implies that $\mu_\bias \ll \mu_0$, analogously to the
independent case treated in Theorem A.
We answer this negatively in Sections \ref{sec:counterexample} and
\ref{sec:bias},
where the
following is proved.

\noindent {\em Notation: } Write $a_n \asymp b_n$ to mean that
there exist positive finite constants $C_1,C_2$ so that $C_1 \leq
a_n/b_n \leq C_2$ for all $n \ge 1$.

\begin{thm} \label{thm:counterex}
  Let $1/2 < \gamma < 1$.  Suppose that the return probabilities $\{u_n\}$
  satisfy $u_n \asymp n^{-\gamma}$ and $\max \{u_i \; : \; i \geq 1\} >
2^{\gamma -1}$.
  \begin{description}
  \item{(i)} If $\theta > {2^{\gamma} \over \max \{u_i \; : \; i\geq 1\}} -
1$, then $\mu_\theta
  \perp \mu_0$.
  \item{(ii)} The bias $\bias$ can be a.s.\ reconstructed from the
  coin tosses $\{X_n\}$, provided $\bias$ is large enough.
  More precisely, we exhibit
  a measurable function $g$ so that,
  for all $\bias > {2^{\gamma} \over \max\{i\; : \; i \ge 1\} } - 1$, we have
  $\bias=g(X)$ $\mu_\bias$-almost surely.
  \end{description}
  \end{thm}
Part (i) is proved, in a stronger form, in Proposition
\ref{prop:cntex}, and
(ii) is contained in Theorem \ref{thm:reconstr} in Section
\ref{sec:bias}, where $g$ is defined.

In Section \ref{sec:counterexample} we provide examples of random walks having
return probabilities satisfying the hypotheses of Theorem
\ref{thm:counterex}.  We provide other examples of Markov chains in
this category in Section \ref{sec:example}.

For this class of examples, Theorem B(ii) and Theorem
\ref{thm:counterex}(i) imply that there is a {\bf phase transition} in
$\theta$: there is a critical $\theta_{c} \in (0,1)$ so that for $\theta <
\theta_c$, the measures $\mu_\theta$ and $\mu_0$ are equivalent,
while for $\theta >
\theta_c$,  $\mu_\theta$ and $\mu_0$ are mutually singular.
  See Section \ref{sec:gc}
for details.  Consequently, there are cases of absolute
continuity, where altering the underlying
Markov chain by introducing delays can
produce singularity.

Most of our current knowledge on the critical parameter
$$\theta_c \; \deq \; \sup\{\theta \, : \,  \mu_\bias \ll \mu_0\}$$ 
is summarized
in the following table.  Choose $r$ such that 
$
u_r \; = \; \max\{u_i \; : \; i\geq 1\}$,
and let
$\theta_s = (\sum_{n=1}^{\infty} u_n^2 )^{-1/2} \wedge
 1$. (The arguments of Harris and Keane~\cite{HK} imply
 that $\theta_s$ is the critical parameter
for $\mu_\theta$ to have a square-integrable density with respect to $\mu_0$.)


\begin{center}
\begin{tabular}{|l|l|}
        \hline
        {\bf asymptotics of $u_n$} & {{\bf critical parameters}} \\
        \hline
         $u_n \asymp n^{-1/2}$ &  $0 = \bias_s= \bias_c$ \\
         $u_n \asymp n^{-\gamma}, \; {1 \over 2} < \gamma < 1$ &
         $0 < \theta_s \leq \theta_c \leq u_r^{-1}2^\gamma - 1$ \\
         $u_n = O(n^{-1})$ & $0 < \theta_s \leq \theta_c 
=
1$ \\
        \hline
\end{tabular}
\end{center}

There are  renewal sequences corresponding to the last row for which
 $0 <\theta_s < \theta_c=1$; see Theorem \ref{thm:nopt} and the remark
following it.

Theorem \ref{thm:counterex}(ii) shows that for certain chains
satisfying $\sum_{n=0}^{\infty}u_n^2 < \infty$,
for $\bias$ large enough, the bias $\theta$ of the coin can be
reconstructed from the observations $X$.  Harris and Keane described
how this can be done for all $\bias$
in the case where $\Gamma$ is the simple random
walk on the integers, and asked whether it is possible whenever
$\sum_n u_n^2 = \infty$.  In Section \ref{sec:le} we answer
affirmatively, and prove the following theorem:

\begin{thm} \label{thm:recon}
  If $\sum_n u_n^2 = \infty$, then there is a measurable function $h$
  so that $\theta = h(X)$ $\mu_\bias$-a.s. for all $\bias$.
\end{thm}
In fact, $h$ is a limit of linear estimators (see the proof given in
Section \ref{sec:le}). Theorem \ref{thm:recon} is extended in
Theorem \ref{thm:recongen}.

%

There are examples of renewal sequences with $\sum_k
u_k^2 < \infty$ which do
not exhibit a phase transition: 

\begin{thm} \label{thm:noptz2}  
If the return probabilities $\{u_n\}$
satisfy 
$
u_n  =  O({n}^{-1})
$,
then $\mu_\bias \ll \mu_0$ for all $0 \leq \bias \leq 1$.
\end{thm}
For example, the return probabilities of (even a delayed) random
walk on $\Z^2$ have $u_k \asymp k^{-1}$.

\noindent
{\em Remark:} The significance of this result is that the asymptotic
conditions on $\{u_n\}$ still holds if the underlying Markov chain is
altered to increase the transition probability from the origin to
itself. 

\noindent
This result is proved in Section \ref{sec:z2}.  It is much
easier to prove that $\mu_\bias$ and $\mu_0$ are always
mutually absolutely continuous in the case where the Markov
chain is ``almost transient'', for example if $u_k \asymp
(k\log k)^{-1}$.  We include the argument for this case
as a warm-up to Theorem \ref{thm:noptz2}.  In particular,
we prove the following theorem:

\begin{thm} \label{thm:nopt} If the return probabilities $\{u_n\}$
satisfy $u_k = O(k^{-1})$, and obey the condition
$$
\sum_{k=0}^n u_k \; = \; o\left(\frac{\log n}{\log\log n}\right) \,,$$
then $\mu_\bias \ll \mu_0$ for all $0 \leq \bias \leq 1$.
\end{thm}

Theorem \ref{thm:nopt} is extended in Theorem \ref{thm:nophase} in Section
\ref{sec:example}, and Proposition \ref{prop:nophaseex} provides examples
of Markov chains satisfying the hypotheses.  Then Theorem
\ref{thm:noptz2} is proved in Section \ref{sec:z2}.

Write $J = \sum_{n=0}^\infty \one_{ \{
\Gamma_n = o\}} \one_{ \{ \Gamma'_n = o\}}$, where
$\Gamma$ and $\Gamma'$ are two independent copies of the
underlying Markov chain.
The key to the proof by Harris and Keane of Theorem B(ii) is the
implication
$$\E[ (1+\theta^2)^J ] < \infty \; \Rightarrow \; \mu_\theta \ll \mu_0 \, .$$
To prove Theorem \ref{thm:noptz2} and Theorem \ref{thm:nopt}
we refine this and show that
$$\E[ (1+\theta^2)^J \; | \; \Gamma ] < \infty \; \Rightarrow \; \mu_\theta
\ll \mu_0 \, . $$

The model discussed here can be generalized by substituting
real-valued random variables for the coin tosses.  We consider the model where
observations are generated with distribution $\measn$ at times when the
chain is away from $o$, and a distribution $\measu$ is used when the
chain visits $o$.

Similar problems of ``random walks on scenery'' were considered by
Benjamini and Kesten in \cite{BK} and by Howard in \cite{Ho1,Ho2}.
Vertices of a graph
are assigned colors, and a viewer, provided only with the sequence of
colors visited by a random walk on the graph, is asked to distinguish
(or reconstruct) the coloring of the graph.

The rest of this paper is organized as follows.  In Section
\ref{sec:defns}, we provide definitions and introduce notation.  In
Section \ref{sec:gc}, we prove a useful general zero-one law, to show
that singularity and absolute continuity of the measures are the only
possibilities.    In Section
\ref{sec:counterexample}, Theorem \ref{thm:counterex}(i) is proved,
while Theorem \ref{thm:counterex}(ii)
is established in Section \ref{sec:bias}.  We prove a more general
version of Theorem \ref{thm:recon} in Section \ref{sec:le}.
In Section \ref{sec:qmgf}, we prove a criterion for
absolute continuity, which is used to prove Theorem \ref{thm:nopt} in
Section \ref{sec:example} and Theorem \ref{thm:noptz2} in Section
\ref{sec:z2}.  A connection to long-range percolation and
some unsolved problems are described in Section \ref{sec:conclusion}.

\section{Definitions.} \label{sec:defns}

Let $\ren = \bin^{\infty}$ be the space of binary sequences.
Denote by  $\renp_{n}$ the $n^{th}$ coordinate projection from $\ren$.
Endow  $\ren$ with the $\sigma$-field  $\HH$ generated by
$\{\renp\}_{n \ge 0}$ and
let $\P$ be a renewal measure on
$(\ren,\HH)$, that is, a measure obeying \be \label{eq:renmeas}
\P[\renp_{0}=1,\renp_{n(1)}=1,\ldots,\renp_{n(m)}=1] =
\prod_{i=1}^{m}u_{n(i)-n(i-1)}, \ee where $u_{n} \deq
\P[\renp_{n}=1].$  We let $\{ T_k \}_{k=1}^{\infty}$ denote the {\em
  inter-arrival times} of the renewal process: If $S_n = \inf\{m >
S_{n-1} : \renp_m = 1\}$ is the time of the $n^{th}$ renewal, then $T_n = S_n
-  S_{n-1}$.  The condition (\ref{eq:renmeas}) implies that
$T_1,T_2,\ldots$ is an \iid sequence. We will use $f_{n}$ to denote
$\P[T_{1} = n]$.

In the introduction we defined $u_n$ as the probability for
a Markov chain $\Gamma$ to return to its initial state at time $n$.
  If $\renp_n =\one_{\{\Gamma_n = o\}}$, then the Markov property
guarantees that
(\ref{eq:renmeas}) is satisfied.  Conversely, any renewal process
$\renp$ can be realized as the indicator of return times of a Markov
chain to its initial state.  (Take, for example, the chain whose value at
epoch $n$ is the time until the next renewal, and consider returns to
$0$.)  Thus we can move freely between these points of view.  For
background on renewal theory, see \cite{F1} or \cite{Ki}.

Suppose that $\measn, \measu$ are two probabilities on
$\R$ which are mutually absolutely continuous,
that is, they share the same null sets.  In
the coin tossing case discussed in the Introduction, these measures
are supported on $\{-1,1\}$.  Given a renewal process,
independently generate observations according to $\measu$ at renewal
times, and according to $\measn$ at all other times.  We describe the
distribution of these observations for various choices of $\measu$.

Let $\R^{\infty}$ denote the space of real sequences,
endowed with the $\sigma$-field $\G$ generated by coordinate projections.
Write $\measu^{\infty}$ for the product probability on $(\R^\infty,\G)$ with
marginal $\measu$.
Let $\Q_{\measu}$ be the measure
$\measn^{\infty}\times\measu^{\infty}\times\P$ on $(\R^\infty \times
\R^\infty \times \ren ,\G \otimes \G \otimes \HH)$.  In the case where
$\eta$ is the coin tossing measure with bias $\bias$, write $\Q_\bias$
for $\Q_\eta$. The random
variables $\np_{n},\sp_{n}$ are defined by
$\np_{n}(\nt,\st,\rent) = \nt_{n}, \; \; \sp_{n}(\nt,\st,\rent) =
\st_{n}.$
Finally, the random variables $X_n$ are
defined by
$$X_{n}=(1-\renp_{n})\np_{n} + \renp_{n}\sp_{n}.$$
The distribution of
$X=\{X_n\}$ on $\R^\infty$ under $\Q_{\measu}$ will be denoted $\muu$.

The natural questions in this setting are: if $\measa$ and $\measb$
are two mutually absolutely continuous measures on $\R$, under what
conditions is $\ma
\perp \mb$?  Under what conditions is $\ma \ll \mb$?  When can $\measu$
be reconstructed from the observations $\{X_n\}$ generated under
$\muu$?  Partial answers are provided in Proposition \ref{prop:cntex},
Theorem \ref{thm:counterex}, Theorem \ref{thm:reconstr}, Theorem
\ref{thm:recongen}, and Theorem \ref{thm:nophase}.

\section{A Zero-One Law and Monotonicity.}  \label{sec:gc}
We use the notation established in the previous section.
Let $\G_n$ be the $\sif$ on $\R^{\infty}$ generated by the first $n$
coordinates.  If $\ma$ and $\mb$ are both restricted to $\G_n$, then
they are mutually absolutely continuous, and we can define the
Radon-Nikodym derivative
$\rho_n = {d\mb \over d\ma}|_{\G_n}$.  Write $\rho$ for $\liminf_{n
  \rightarrow \infty}\rho_n$; the Lebesgue Decomposition Theorem (see
  Theorem 4.3.3 in \cite{Du}) implies
that for any $A\in \G$, \be
\label{eq:lebdec}
\mb[A] = \int_A \rho d\ma + \mb^{sing}(A)
= \int_A \rho d\ma + \mb[\{\rho = \infty\} \cap A],\ee
where $\mb^{sing} \perp \mu_{\measa}$. Thus
to prove that $\mb \ll \ma$, it is enough to show that \be
\label{eq:condrho} 1= \mb[x: \rho(x) < \infty ] = \Q_{\measb}[\rho(X)
< \infty].  
\ee 
For any process $\Gamma$, let $\shift_n \Gamma =
(\Gamma_{n},\Gamma_{n+1},\ldots)$, and let $\tail(\Gamma) =
\bigcap_{n=1}^{\infty}\sigma(\shift_n \Gamma)$ be the tail $\sif$.
\begin{lem}[Zero--One Law] \label{lem:zo} The tail $\sigma$-field
  $\tail(\np,\sp,\renp)$, and hence $\tail(X)$, is
  $\Q_{\measu}$-trivial.
That is, $A \in \tail(\np,\sp,\renp)$ implies
  $\Q_{\measu}(A) \in \{0,1\}$.
\end{lem}
\proof{Proof}  By the Kolmogorov Zero-One Law, $\tail(\np)$ and
$\tail(\sp)$ are trivial.  The inter-arrival times $\{T_n\}$ form an
\iid sequence, and clearly $\tail(\renp) \subset {\cal
  E}(T_1,T_2,\ldots)$, where ${\cal E}$ is the exchangeable $\sigma$-field.
  The Hewitt-Savage Zero-One law implies that
${\cal E}$, and hence $\tail(\renp)$, is trivial.

Let $f$ be a bounded $\tail(\np,\sp,\renp)$-measurable function on
$\R^\infty\times
\R^\infty \times \ren$ which can be written as
\be \label{eq:funcform}
f(\nt,\st,\rent) = f_1(\nt)f_2(\st)f_3(\rent)\, . \ee
By independence of $\np$,$\sp$, and $\renp$, and triviality of
$\tail(\np)$,$\tail(\sp)$, and $\tail(\renp)$, it follows that
$$
\E [f_1(\np)f_2(\sp)f_3(\renp)] = \E f_1(\np) \E f_2(\sp) \E f_3(\renp)
= f_1(\np) f_2(\sp) f_3(\renp) \; a.s.
$$
  Consequently, for all functions
of the form (\ref{eq:funcform}),
\be \E f(\np,\sp,\renp) = f(\np,\sp,\renp) \; a.s. \label{eq:funcres}
\ee
The set of bounded functions of the form (\ref{eq:funcform}) is closed
under multiplication, includes
the indicator functions of rectangles $A \times B \times C$ for $A,B \in \HH$
and $C \in \G$, and these rectangles generate the $\sif$
$\G\times\G\times\HH$.  Since the collection of
bounded functions satisfying (\ref{eq:funcres}) form a
monotone vector space, a Monotone Class Theorem implies that all
bounded $\G\times\G\times\HH$-measurable functions obey
(\ref{eq:funcres}).
We conclude that $\tail(\np,\sp,\renp)$ is
trivial. 
\qed
\endproof

\begin{prop}\label{prop:eitheror} Either $\mb$ and $\ma$ are mutually
absolutely continuous, or $\mb \perp \ma$.
\end{prop}
{\sc Proof.}  Suppose that $\mu_\measb \not \perp \mu_\measa$.
  {}From (\ref{eq:lebdec}), it must be that $\rho < \infty$ with
positive $\mu_\measb$ probability.
Because the event $\{\rho < \infty\}$ is in $\tail$,
Lemma \ref{lem:zo} implies
$\rho < \infty$ $\mb$-almost surely.  Using (\ref{eq:lebdec}) again,
we have that $\mu_\measb \ll \mu_\measa$.  The same argument
with the roles of $\measa$ and $\measb$ reversed, yields that
$\mu_\measa \ll \mu_\measb$ also.
%
%
%
%
\qed

We return to the special case of coin tossing here, and justify our
remarks in the introduction that for certain sequences $\{u_{n}\}$,
there is a phase transition.  In particular, we need the following
monotonicity result.
\begin{prop} \label{prop:mono} Let $\theta_1 < \theta_2$.  If
$\mu_{\theta_1} \perp
  \mu_{0}$, then $\mu_{\theta_2} \perp \mu_0$.
\end{prop}
{\sc Proof.}  Couple together the processes $X$ for all $\bias$: At each
epoch $n$,
generate a variable $V_n$, uniformly distributed on $[0,1)$.  
If $\renp$ is a renewal process independent of $\{V_n\}$,
define $X^\bias$ by
\be \label{eq:couple} X^\bias_n = \left\{ \begin{array}{ll}
+1 & \mbox{if } V_n \leq {1+\bias\renp_n \over 2}
\vspace{0.05in}\\
-1 & \mbox{if } V_n > {1+\bias\renp_n \over 2} \end{array} \right. \ee
Then $X^{\bias_1} \leq X^{\bias_2}$ for $\bias_1 < \bias_2$, and
$X^{\bias}$ has law $\mu_\bias$ for all $\bias \in [0,1]$.
Thus $\mu_{\theta_{2}}$
stochastically dominates $\mu_{\theta_{1}}$.

Suppose now that $\mu_{\theta_{1}} \perp \mu_{0}$.  Then
(\ref{eq:lebdec}) implies that \be \mu_{\theta_{1}}[\rho_{\theta_{1}}
= \infty] = 1 \; \mbox{and } \; \mu_{0}[\rho_{\theta_{1}} = 0] =1.
\label{eq:twoconds} \ee Because the functions
$$\rho_n(x) = \int_\Upsilon \prod_{k=0}^n (1+\bias x_k
\renp_k)d\P(\renp)$$
are increasing in $x$, it follows that $\rho$ is an increasing
function and the event $\{\rho = \infty\}$ is an increasing event.
Because $\mu_{\theta_2}$ stochastically dominates
$\mu_{\theta_1}$, we have
\be\mu_{\theta_{2}}[\rho_{\theta_{1}} = \infty] = 1.
\label{eq:onecond} \ee
Putting together (\ref{eq:onecond}) and the second part of
(\ref{eq:twoconds}) shows that we have decomposed $\R^\infty$ into the
two disjoint sets $\{\rho_{\theta_1} = 0\}$ and
$\{\rho_{\theta_1} = \infty \}$ which satisfy
$$\mu_{0}[\rho_{\theta_1} = 0] =1 \quad \mbox{and} \quad
\mu_{\theta_2}[\rho_{\theta_1} = \infty ] =1\,.$$
In other words, $\mu_{\theta_{2}} \perp \mu_{0}$. \qed

Consequently, it makes sense to define for a given renewal sequence
$\{u_{n}\}$ the {\em critical bias} $\bias_{c}$ by $$\theta_{c} \deq
\sup\{\bias\leq 1: \mu_{\theta} \ll \mu_{0}\}.$$
We say there is a {\em
  phase transition} if $0 < \bias_{c} < 1$.  The results of
Harris and Keane say $\sum u_{n}^{2} = \infty$ implies $\theta_{c}=1$
and there is no phase transition.  In Section
\ref{sec:counterexample}, we provide examples of
$\{u_{n}\}$ with $\sum_n u_n^2 < \infty$ having a phase transition.
In Section \ref{sec:example}, we provide examples
with $\sum_n u_n^2<\infty$ without a phase transition.

\section{Existence of Phase Transition.} \label{sec:counterexample}
In this section, we confine our attention to the coin tossing
situation discussed in the Introduction.  In this case, $\measn$ and $\measa$
are both the probability on $\coins$ with zero mean,
and $\measb$ is the probability with
mean $\bias$ (the $\bias$-biased coin).  The distributions
$\mu_\measa$ and $\mu_\measb$ are denoted by $\mu_0$ and $\mu_\bias$
respectively. Let $U_n \deq \sum_{k=0}^n u_k$.


\begin{prop} \label{prop:cntex}
Let $\{u_n\}$ be a renewal sequence with 
$$
\sum_{k=0}^n u_k \; = \; U_n
\; \asymp \; n^{1-\gamma}\ell(n)\,,$$
for
${1\over 2} < \gamma < 1$ and $\ell$ a slowly varying function.
If  
$$(1+\theta)\max \{u_i \; : \; i \geq 1\} \; > \; 2^{\gamma} \,,$$
then $\mu_{\theta} \perp \mu_0$.
\end{prop}
{\bf Remark. } The conditions on $\bias$ specified in the statement
above are not vacuous.  That is, there are examples where the
lower bound on $\bias$ is less than $1$.  There are random walks with
return times obeying $u_n \asymp n^{-\gamma}$, as shown in Theorem
\ref{prop:stable}.  By introducing delays at the origin, $u_1$ can be
made to be close to $1$, so that $2u_1 > 2^\gamma$.

\noindent
{\sc Proof.}
Let $\E$ denote expectation with respect to the renewal measure $\P$
and let $\E_\theta$ denote expectation with respect to $\Q_\bias$.
Let 
$u_r = \max \{u_i \; : \; i \geq 1\}$
and assume for now that $r=1$.
Let 
$b = \frac{1}{2}(1+\theta)$ 
and 
$k(n) = \lfloor (1 + \ep)\log_2n \rfloor$,
where $\ep$ is small enough that 
$(1+\ep)(-\log_2 u_1 b) <1 - \gamma$.
Define $A_j^n$ as the event that at all times
$i \, \in \, \left[ jk(n), (j+1)k(n) \right)$
there are renewals and the coin lands ``heads'', i.e.,
$$
A_j^n \, \deq \, \bigcap_{\ell=0}^{k(n)-1}\{\renp_{jk(n)+\ell} = 1 \;
\mbox{ and }
X_{jk(n)+\ell} = 1\}.
$$
Let 
$D_n \deq \sum_{j=1}^{n/k(n)} \one_{A_j^n}$,
and 
$$
c(n) \; \deq \; \Q_\bias [A_{j}^{n}|\renp_{jk(n)}=1] \; = \;
b (u_1 b)^{k(n)-1} = u_1^{-1}(u_1 b)^{k(n)} \,.$$
Note that we have defined things so that
$c(n) \asymp n^{-p}$, where $p < 1-\gamma$.  Then
\be \label{eq:ex1}
\E_\bias D_n = \sum_{j=1}^{n/k(n)} u_{jk(n)} b (u_1 b)^{k(n)-1} = 
c(n) \sum_{j=1}^{n/k(n)} u_{jk(n)} .\ee
We need the following simple lemma:
\begin{lem} \label{lem:lemren} 
    For all $r \geq 0$,
\be
u_r + u_{r+k} + \cdots + u_{r+mk} \leq u_0 + u_k + \cdots + u_{mk} \,.
\label{eq:renlemma} \ee
\end{lem}
{\sc Proof.} Recall that $u_0 = 1$.
Let $\tau^* = \inf\{j \geq 0:
\renp_{r+jk}=1\}$.
Then 
\bean
\E[ \sum_{j=0}^m \renp_{jk+r} | \tau^* ]
& = & (1+ u_{k} + \cdots + u_{(m-\tau^*)k})\one_{\{\tau^* \leq m\}}\\
& \leq & u_0 + u_k +\dots + u_{mk} \, .
\eean
Taking expectation proves the lemma. \qed

By this lemma,
\be
\sum_{j=0}^{n/k(n)} u_{jk(n)} \geq {1\over k(n)}
\sum_{j=0}^n u_j  = {U_n \over k(n)} \, ,
\nonumber \ee
and thus
\be 
\sum_{j=1}^{n/k(n)} u_{jk(n)} \geq { U_n \over k(n)} - 1 \asymp
 {U_n \over k(n)} \asymp n^{1-\gamma}{\ell(n)
\over k(n)} \,.
\label{eq:ex2} \ee
Combining (\ref{eq:ex1}) and (\ref{eq:ex2}), we find that
$$
\E_\bias D_n \geq C_1 n^{-p} n^{1-\gamma}{\ell(n) \over k(n)} =
C_1{n^{1-\gamma-p}}{\ell(n) \over
k(n)}\,.
$$

Since $1-\gamma-p>0$, it follows that $\E_\bias D_n \rightarrow \infty$.

Also,
\bea
\E_\bias D_{n}^{2} & = &
\sum_{i=1}^{n/k(n)}\Q_\bias[A_{i}^{n}] +
{2}\sum_{i=1}^{n/k(n)}\sum_{j=i+1}^{n/k(n)}\Q_\bias [A_{j}^{n}|A_{i}^{n}]
\Q_\bias [A_{i}^{n}] \nonumber \\
& = & \E_\bias D_{n} + {2}
\sum_{i=1}^{n/k(n)}\sum_{j=i+1}^{n/k(n)} c(n) u_{k(n)(j-i-1)+1}
c(n) u_{k(n)i}
\nonumber \\
& \leq & \E_\bias D_{n} + {2}c(n)^2
\sum_{i=1}^{n/k(n)} u_{k(n)i} \sum_{j=0}^{n/k(n)} u_{k(n)j+1}
\nonumber \\
& \leq & \E_\bias D_{n} + {2}c(n)^2
\sum_{i=1}^{n/k(n)} u_{k(n)i} \sum_{j=0}^{n/k(n)} u_{k(n)j}
\label{eq:uselem} \\
& \leq & \E_\bias D_{n} + {2}c(n)^2
\sum_{i=1}^{n/k(n)} u_{k(n)i} \sum_{j=1}^{n/k(n)} u_{k(n)j}
+ 2u_0 c(n) \E_\bias D_n 
\label{eq:termj0} \\
& \leq & \E_\bias D_{n} +
{2 c(n)^{2}} \left(\sum_{i=1}^{n/k(n)}u_{k(n)i}\right)^{2}
+ 2u_0 c(n) \E_\bias D_n 
\nonumber \\
& \leq & C(\E_\bias D_{n})^{2} \label{eq:smm} 
\eea
(\ref{eq:uselem}) follows from Lemma \ref{lem:lemren}, and
the last term in (\ref{eq:termj0}) comes from the contributions when 
$j=0$.

If $A_n$ is the event that there is a run of length $k(n)$ after epoch
$k(n)$ and before $n$, then (\ref{eq:smm}) and the second moment
inequality yield
$$\Q_\bias[A_{n}] \geq \Q_\bias[D_{n} >0 ] \geq \frac{(\E_\bias
  D_{n})^{2}}
{\E_\bias D_{n}^{2}} \geq \frac{1}{C} > 0.$$
Finally, we have
$$\Q_\bias[\limsup A_n] \geq \limsup \Q_\bias[A_n] > 0,$$
and by the Zero-One Law
(Lemma \ref{lem:zo}) we have that $\Q_\bias[\limsup A_n] = 1$.  A theorem
of Erd\H{o}s and R\'enyi (see, for example, Theorem 7.1 in \cite{Rev})
states that under the measure $\mu_0$, $L_n/\log_2n \rightarrow 1$,
where $L_n$ is the length of the longest run before epoch $n$.  But
under the measure $\mu_\bias$, we have just seen that we are
guaranteed to, infinitely often, see a run of length
$(1+\epsilon)\log_2 n$ before time $n$.

If $u_1 \neq \max \{u_i \; : \; i \geq 1\}$,  consider the renewal process
$\{\renp_{nr}\}_{n=0}^{\infty}$ and the sequence
$\{X_{nr}\}_{n=0}^{\infty}$, where $u_r = \max \{u_i \; : \; i \geq 1\}$.
Apply the
proceeding argument to this subsequence to distinguish between
$\mu_\bias$ and $\mu_0$. \qed

\begin{prop} \label{prop:stable}
There exists a renewal measure $\P$ with $u_{n} \sim
  Cn^{-\gamma}$ for $1/2 < \gamma < 1$.
\end{prop}
{\sc Proof.}  For a distribution function $F$ to be in the
domain of attraction of a stable law, only the asymptotic behavior
of the tails $F(t), 1- F(-t)$ is relevant (see, for
example, Theorem 8.3.1 in \cite{BGT}).  Thus if the symmetric stable
law with exponent $1/\gamma$ is discretized so that it is supported on
$\Z$, then the modified law $F$ is in the domain of attraction of this
stable law.  Then if $\Gamma$ is the random walk with increments
distributed according to $F$, Gnedenko's Local Limit Theorem (see
Theorem 8.4.1. of \cite{BGT}) implies that
$$\lim_{n\rightarrow\infty}|n^{\gamma}\P[\Gamma_{n} =0] -g(0)| = 0,$$
where $g$ is the density of the stable law.  Thus if $\renp_{n}\deq
\one_{\{\Gamma_{n}=0\}}$, then $\{\renp_{n}\}$ form a renewal sequence
with $u_{n}\sim Cn^{-\gamma}.$ \qed

For a sequence to satisfy the hypotheses of Proposition
\ref{prop:cntex} and \ref{thm:counterex}, we also need that
$\max \{u_i \; : \; i \geq 1\} > 2^{\gamma-1}$.
By
introducing a delay at the origin for the random walk $\Gamma$ in
Proposition \ref{prop:stable}, $u_{1}$ can be made
arbitrarily close to $1$.  Thus there do exist Markov chains which
have $0 < \bias_c < 1$.

An example of a Markov chain with $U_n \asymp n^{1/4}$ will be
constructed by another method in Section \ref{sec:example}.

\section{Determining the bias $\bias$.} \label{sec:bias}
In this section we refine the results of the previous section and
give conditions that allow reconstruction of the bias from the
observations.

For $a\geq 1$, let
\be \label{eq:lambdadef} 
\legt(a) \; \deq \; \lim_{m \rightarrow \infty} 
\frac{ -\log_2 \P[ T_1 + \cdots + T_m \leq ma]}{m} \,. 
\ee
($\legt(a) = \infty$ for $a<1$ (since each $T_i \geq 1$), hence we 
restrict attention to when $a \geq 1$.)

Because $\E T_i=\infty$, Cram\'er's Theorem
(see, e.g., \cite{DeZ})
implies that $\legt(a) > 0$ for all $a$.
Since 
$\lim_{a \uparrow \infty}\P[T_1 \leq a]
= 1$, 
it follows that
$\lim_{a \uparrow \infty}\legt(a)
= 0$.  
Also, $\legt(1) = - \log_2 u_1$.  

It is convenient to reparameterize so that
we keep track of 
$\param \deq \log_2(1+\bias)$ instead of $\bias$
itself. Let 
\be \label{eq:defnpsi}
\widehat{\psi}(\param,\xi) \deq \xi\cdot(\param -\legt(\xi^{-1})) 
\; \mbox{and} \; \psi(\param) \deq \sup_{0 < \xi \leq 1}
\widehat{\psi}(\param,\xi).
\ee
Observe that $\lim_{\xi \rightarrow 0}\widehat{\psi}(\param,\xi) =  0$.
For $\ep>0$ small enough so that $\legt(\ep^{-1}) < \frac{\param}{2}$, 
$$
\widehat{\psi}(\param,\ep) > \ep (\param - \frac{\param}{2})
= \ep \frac{\param}{2} > 0\,.
$$
Hence, 
The maximum of $\widehat{\psi}(\param,\cdot)$ over
$(0,1]$ is attained, so we can
define
$$
\xi_0 (\param) \; \deq \; 
\inf\{0 < \xi \leq 1 \; : \; \widehat{\psi}(\param,\xi) = 
\psi(\param)\} \,. 
$$
We show now that
$\widehat{\psi}(\param,\xi_0) > \widehat{\psi}(\param,1)$, 
a fact which we will use later (see the
remarks following Theorem \ref{thm:reconstr}).
Let $\ell = \min\{n>1:f_n>0\}$, and note that $f_1 = u_1$.
If in the interval $[0,k(1+\lfloor \ep\ell \rfloor )\,]$ there are
$k-\lfloor \ep k \rfloor$
inter-renewal times of length $1$ and $\lfloor \ep k \rfloor$
 inter-renewal times of
length $\ell$, then in particular there are at least $k$ renewals.
Consequently,
\be \P[T_1 + \cdots + T_k \leq k(1+\ep \ell)] \; \geq \;
{k \choose \lfloor \ep k\rfloor}
 f_1^{k-\lfloor \ep k \rfloor} f_\ell^{\lfloor \ep k
\rfloor } \label{eq:ent1} \ee
Taking logs, normalizing by $k$, and then letting $k \rightarrow
\infty$ yields
\bean
-\legt(1+\ep\ell) & = & \lim_{k \rightarrow \infty}
k^{-1}\log_2 \P[T_1 + \cdots + T_k \leq k(1+\ep \ell)] \\
& \geq & h_2(\ep) + \log_2 f_1 + \ep \log_2(f_\ell/f_1) \,,
\eean
where $h_2(\ep) = \ep \log_2 \ep^{-1} + (1 - \ep)\log_2(1-\ep)^{-1}$.
Therefore
\bea
\psi(\param,{1 \over 1+\ep\ell}) - \psi(\param,1) & = &
 {1 \over 1+\ep\ell}\param - {1 \over 1+\ep\ell}
 \legt(1 + \ep\ell) - \param - \log_2
 f_1 \label{eq:diff} \\
 & \geq & {1 \over 1+\ep\ell}\left\{ -\ep(\ell\param +
 \log_2(f_\ell/f_1))
  + h_2(\ep) \right\} \label{eq:diff2}
\eea
Thus for $\ep$ bounded above, the
left-hand side of (\ref{eq:diff}) is bounded below by
$C_1 ( h_2(\ep) - C_2\ep)$.  
Since the derivative of $h_2$ tends to infinity near $0$, there is
a positive $\ep$ where the difference is strictly positive.  Thus, the
maximum of $\whpsi(\param,\cdot)$ is {\em not} attained at $\xi = 1$.

Finally, $\psi$ is
strictly increasing: let $\param < \param'$, and observe that
$$\psi(\param') = \whpsi(\param',\xi_0(\param')) \geq
\whpsi(\param',\xi_0(\param)) > \whpsi(\param,\xi_0(\param)) =
\psi(\param).$$
\begin{thm}\label{thm:reconstr}
Recall that 
$$
\P[X_k = 1 \; | \; \renp_k = 1] = 2^{-1}(1+\bias) =
2^{\param - 1}, \, \mbox{for } \param \deq \log_2
  (1+\theta) \,.
$$
Let 
$$
R_n = \sup\{m:X_{n+1} = \cdots =
X_{m+n} = 1\} \quad \mbox{and} \quad
\ls(X) = \limsup_n R_n(\log_2 n)^{-1} \,.
$$
Suppose that ${1\over 2} < \gamma < 1$ and $\ell$ is a slowly varying
function. If
$U_n \asymp n^{1-\gamma}\ell(n)$, then
$$\ls(X) = \frac{1 - \gamma}{1 - \psi(\param)} \bigvee 1,$$
where $\psi$ is the strictly monotone function defined in (\ref{eq:defnpsi}).

In particular, for $\param > \psi^{-1}(\gamma)$
(equivalently, $\bias \geq 2^{\psi^{-1}(\gamma)} - 1$),
we can recover $\param$ (and hence $\bias$) from $X$:
$$
\param  \; = \; \psi^{-1}\left(1 -
{1-\gamma \over \ls(X)}\right) \,.
$$
\end{thm}
{\bf Remark.}  Suppose $u_1= \max \{u_i \; : \; i \geq 1\}$.  Since
$\psi(\param) >
\widehat{\psi}(\param,1)$, 
(see the comments before the statement of Theorem \ref{thm:reconstr})
we have that
\be \label{eq:str}
\psi(\param)  \; > \; \param + \log_2 u_1 \,.
\ee
Substituting 
$\psi^{-1}(\gamma)$ for $\param$ in (\ref{eq:str}) yields 
$$\psi^{-1}(\gamma) \; < \; \gamma - \log_2 u_1 \,.$$
Thus 
\be \label{eq:comp}
2^{\psi^{-1}(\gamma)} - 1 \; < \; 2^{\gamma -
  \log_2u_1} -1.
\ee 
The right-hand side of (\ref{eq:comp}) is the
upper bound on $\bias_c$ obtained in Proposition \ref{prop:cntex},
while the left-hand side is the upper bound given by Theorem
\ref{thm:reconstr}.  Thus this section strictly improves the results
achieved in the previous section.

\noindent
{\sc Proof.}
Let $\c = (1-\gamma)/(1-\psi(\param))$.
We begin by proving that $\ls(X) \leq \c \vee 1$, or equivalently,
that \be \label{eq:part1}
\forall c> \c\vee 1 , \; \Q_\bias[R_n \geq c\log_2 n \; i.\ o.\ ] = 0. \ee
Fix $c > \c \vee 1$.  If $k(n,c) = k(n) \deq \lfloor c\log_2 n \rfloor$,
then it is enough to
show that \be
\Q_\bias[ \limsup_n \{X_{n+1} = \cdots = X_{n+k(n)} =1 \} ] = 0.
\label{eq:bnd1}\ee

Let $E_n$ be the event $\{X_{n+1} =\cdots = X_{n+k(n)} =1 \}$, and
define 
$$F_n \; \deq \; 
\inf\{m>0: \renp_{n+m} =1\}$$ as the waiting time at $n$
until the next renewal (the residual lifetime at $n$).
We
have \be \Q_\bias[E_n] \leq \Q_\bias
[E_n \; | \; F_n > k(n) ] + \sum_{m=1}^{k(n)} \Q_\bias[E_n
\; | \; F_n = m]\Q_\bias[F_n = m] . \label{eq:condfn} \ee Notice that  $$
\{F_n = m \} = \{\renp_{n+1} = \cdots = \renp_{n+m-1}=0, \;
\renp_{n+m}=1\},$$
and consequently we have 
\be
\begin{array}{c}
    \Q_\bias[E_n \; | F_n = m, \; \renp_{n+m+1},\ldots,\renp_{n+k(n)} ] \\
    \; = \; 2^{-k(n)}(1+\bias)^{1+\renp_{n+m+1}+\cdots+\renp_{n+k(n)}} 
\end{array}\quad .
\label{eq:encd} 
\ee
Taking expectations over $(\renp_{n+m+1},\ldots,\renp_{n+k(n)})$ in
(\ref{eq:encd}) gives that \bea \Q_\bias[E_n \; | \; F_n = m] & = &
2^{-k(n)}\E[ (1+\bias)^{1+\renp_{n+m+1}+\cdots+\renp_{n+k(n)}}
\; | \; \renp_{n+m}=1 ] \nonumber \\
& = & 2^{-k(n)}\E
[(1+\bias)^{1+\renp_{1}+\cdots+\renp_{k(n)-m}}]\label{eq:strmk} \eea
The equality in (\ref{eq:strmk}) follows from the renewal property,
and clearly the right-hand side of (\ref{eq:strmk}) is
maximized when $m=1$.   
Therefore the right-hand side
of (\ref{eq:condfn}) is bounded above by 
\be 
2^{-k(n)} + (U_{n+k(n)}-U_n) \Q_\bias[E_n \; | \; \renp_{n+1} = 1].
\label{eq:fb}
\ee 
We now examine the probability $\Q_\bias[E_n \; | \; \renp_{n+1} = 1]$
appearing on the right-hand side of
(\ref{eq:fb}).  Let 
$N[i,j] \deq \sum_{k=i}^j \renp_k$ be
the number of renewals appearing between times $i$ and $j$.  In the
following, $N=N[n+1,n+k(n)]$.  We have
\bea 
\Q_\bias[E_n \; | \; \renp_{n+1} = 1] 
& = & 2^{-k(n)} \E[ (1+\bias)^N \given \renp_{n+1} = 1 ] 
\nonumber \\
& = & 
\E[2^{k(n)(-1 + \param N/k(n))} \given \renp_{n+1}=1] \,.
\label{eq:cp1}
\eea
By conditioning on the possible values of $N$, (\ref{eq:cp1}) is
bounded by 
\be \label{eq:poa1} 
\sum_{m=1}^{k(n)} 2^{k(n)(-1 + \param m/k(n))} 
\P [ T_1 + \cdots T_m \leq k(n) ] \,. 
\ee
By the superadditivity of $\log\P[T_1+\cdots+T_m \leq m a]$,
the
probabilities in the sum in (\ref{eq:poa1}) are bounded above by
$2^{-m\legt(k(n)/m)}$.  Consequently, (\ref{eq:poa1}) is dominated by
\bean
\sum_{m=1}^{k(n)} 2^{k(n)(-1 + m/k(n)(\param - \legt(k(n)/m)))}
& \leq & \sum_{m=1}^{k(n)} 2^{k(n)(-1 + \whpsi(\param,m/k(n)))} 
\nonumber \\
& \leq & k(n) 2^{k(n)(\psi(\param)-1) }
\eean
Hence, returning to (\ref{eq:fb}), \bea
\nonumber
\Q_\bias[E_n]
& \leq & 2^{-k(n)} + (U_{n+k(n)}-U_n) k(n) 2^{-k(n)(1-\psi(\param))}
\\ & \leq & 2n^{-c} + 2k(n)(U_{n+k(n)}-U_n) n^{-c(1-\psi(\param))}.
\label{eq:sum1} \eea
Let $q = c(1-\psi(\param))$, and since $c > \zeta \vee 1$, we have that
$q+\gamma > 1$.   Letting $m(n)=n+k(n)$, since $m(n) \geq n$, we have
\bea
\sum_{n=1}^L k(n)U_{n+k(n)}n^{-q}  &\leq &
\sum_{n=1}^L k(m(n))U_{m(n)}(m(n) - k(n))^{-q}\nonumber \\
&\leq &\sum_{n=1}^L k(m(n))U_{m(n)}(m(n) - k(m(n)))^{-q}\nonumber \\
&\leq &\sum_{m=1}^{L + k(L)}k(m)U_m(m-k(m))^{-q} \, .\label{eq:chvar}
\eea
Then, using (\ref{eq:chvar}), it follows that
\bea
\sum_{n=1}^L k(n)(U_{n+k(n)}-U_n)n^{-q}
& \leq & \sum_{n=1}^Lk(n)U_n\left((n-k(n))^{-q} - n^{-q}\right) 
\nonumber \\
&& + \sum_{n=L+1}^{L+k(L)} k(n)U_n (n-k(n))^{-q} \,. 
\label{eq:cv1}
\eea
Since $a^{-q} - b^{-q} \leq C (b-a) a^{-1-q}$, and $U_n \leq
Cn^{1-\gamma}$, the
right-hand side of (\ref{eq:cv1}) is bounded above by
\be
\begin{array}{l}
C_1 \sum_{n=1}^L k(n)n^{1-\gamma} k(n) (n-k(n))^{-q -1} \\
+ \; C_2 k(L)k(L+k(L)) (L+k(L))^{1-\gamma}(L-k(L))^{-q} 
\end{array} \quad .
\label{eq:aip} \ee
We have that (\ref{eq:aip}), and hence (\ref{eq:cv1}), is bounded 
above by
\be
C_3 \sum_{n=1}^L k(n)^2 n^{-(q+\gamma)} + o(1) \,.
\label{eq:sumbyprts}
\ee
Since $q+\gamma > 1$, (\ref{eq:sumbyprts}) is bounded as $L
\rightarrow \infty.$
We conclude
that (\ref{eq:sum1}) is summable.
Applying the Borel-Cantelli lemma
establishes (\ref{eq:bnd1}).

We now prove the lower bound, $\ls(X) \geq \c \vee 1$.

It is convenient to couple together monotonically the processes
$X^\theta$
for different $\theta$.  See (\ref{eq:couple})
in the proof of Proposition \ref{prop:mono}
for the construction of the coupling, and let $\{V_i\}$ be the \iid
uniform random variables used in the construction.

First, using the coupling, we have that
$\ls(X^\bias) \geq \ls(X^0)=1$.  Hence, 
$$
\mu_\bias[x:\ls(x)\geq 1] \; = \; 1 \,.
$$

It is enough to show
that if $c < \c$, then 
$$
\Q_\bias[R_n \geq k(c,n) \; i.\ o.\ ] = 1 \,.
$$
Fix $\param$, and
write $\xi_0$ for $\xi_0(\param)$.

Let $\tau_i = \tau_i^n$ be the time of the $\lfloor
\xi_0k(n)\rfloor^{th}$ renewal
after time $ik(n)-1$.  The event $G_i^n$ of a {\em good run} in the
block $I^n_i = [ik(n),(i+1)k(n)-1]\cap \Z^+$ occurs when
\begin{enumerate}
\item there is a renewal at time $ik(n)$: $\renp_{ik(n)} =1$,
\item there are at least $\xi_0k(n)$ renewals in $I_i$: $\tau_i \leq
  (i+1)k(n)-1$,
\item until time $\tau_i$, all observations are ``heads'': $X_j=1 \,
\mbox{for }
  ik(n) \leq j \leq \tau_i$,
\item $V_j \leq 1/2 \; \mbox{for } \tau_i < j \leq (i+1)k(n)-1$.
\end{enumerate}
The importance of the coupling and the last condition is that a good
run in $I_i$ implies an observed run ($X_j = 1 \, \forall \, j \in
I_i$).

Let $N_i=N[I_i]$.  The probability of $G^n_i$ is given by \be \Q_\bias
[G^n_i] =
2^{-k(n)}(1+\bias)^{\xi_0 k(n)} p_i u_{ik(n)},
\label{eq:pgn} \ee
where $p_i \deq \P[N_i \geq \xi_0 k(n) \given \renp_{ik(n)} =1]$
is the probability of at least $\xi_0 k(n)$ renewals in the interval
$I_i$, given that there is a renewal at $ik(n)$.  Note that $p_i
\equiv p_1$ for all $i$, by the renewal property.

Following the proof of Proposition \ref{prop:cntex}, we define $D_n =
\sum_{j=1}^{n/k(n)} \one_{G_j^n}$, and compute the first and second
moments of $D_n$.  Using (\ref{eq:pgn}) gives
\be
\E_\bias[ D_n ] =
2^{-k(n)}(1+\bias)^{\xi_0 k(n)} p_1 \sum_{j=1}^{n/k(n)} u_{jk(n)}.
\label{eq:edn1} \ee
Since $c < \zeta = \frac{1-\gamma}{1-\psi(\param)}$, we also have for
some $\ep > 0$ that
\be
c < \frac{1 - \gamma}{1 + \ep\xi_0 - \psi(\param)} \,.
\label{eq:someroom}
\ee
By definition of $\legt$,
we can bound below the
probability $p_1$: For $n$ sufficiently large, 
\bea 
p_1 & = & \P[ N_1 \geq \xi_0 k(n) \given \renp_{k(n)} =1] 
\nonumber \\
& = & \P[T_1 + \cdots + T_{\xi_0 k(n)} \leq k(n)] 
\nonumber \\
& \geq & 2^{-\xi_0 k(n)(\legt(\xi_0^{-1}) + \ep)} \,,
\label{eq:ld1} 
\eea
where $\ep > 0$ is arbitrary.  Thus, plugging (\ref{eq:ld1}) into
(\ref{eq:edn1}) shows that for $n$ sufficiently large,
\bea
\E_\bias[ D_n ] &\geq&
2^{-k(n)}(1+\bias)^{\xi_0 k(n)} 2^{-\xi_0k(n)(\legt(\xi_0^{-1}) + \ep)}
\sum_{j=1}^{n/k(n)} u_{jk(n)}
\nonumber\\
&=& 2^{k(n)(-1 - \ep\xi_0 + \param\xi_0 - \xi_0\legt(\xi_0^{-1}))}
\sum_{j=1}^{n/k(n)} u_{jk(n)} \nonumber \\
& \geq &  2^{-1} n^{-q}
\sum_{j=1}^{n/k(n)} u_{jk(n)} \; ,
\label{eq:almodone} \eea
where $q=(1 + \ep\xi_0 -\psi(\param))c$.  By 
(\ref{eq:someroom}),
$1 -\gamma -q > 0$.
Using (\ref{eq:ex2}), $\sum_{j=1}^{n/k(n)} u_{jk(n)} \asymp
{ \ell(n) \over k(n)} n^{1-\gamma}$, in
(\ref{eq:almodone}), gives that for $n$ large enough,
$$\E_\bias[D_n] \geq C_3 {\ell(n)\over k(n)}n^{1-\gamma - q}
\stackrel{n \rightarrow \infty}{\longrightarrow} \infty.$$

We turn now to the second moment, which we show is bounded by a
multiple of the square of the first moment.
\be \label{eq:fpsm}
\E_\bias [D_n^2 ]
= 2\sum_{i=1}^{n/k(n)} \sum_{j=i+1}^{n/k(n)} \Q_\bias [G_i^n \cap
G^n_j] + \E_\bias[D_n].
\ee
We compute the probabilities appearing in the
sum by first conditioning on the renewal process:
\bea
\Q_\bias[G_i^n \cap G^n_j \; | \; \renp] 
& = & \left( 2^{-k(n)}(1+\bias)^{\xi_0 k(n)}
\right)^2 \one_{\{ N_i \geq \xi_0 k(n), \; \renp_{ik(n)}=1 \}}
\nonumber \\
&& \times \; \one_{\{ N_j \geq \xi_0 k(n), \; \renp_{jk(n)}=1 \}} \,.
\label{eq:intcp} 
\eea
Taking expectations of (\ref{eq:intcp}), if
$d(n)=2^{-k(n)}(1+\theta)^{\xi_0 k(n)}$, then
\bea
\Q_\bias[G_i^n \cap G_j^n] & = & d(n)^2 \P[N_j \geq \xi_0 k(n),
\renp_{jk(n)}=1
\; \mbox{and} \; N_i \geq \xi_0 k(n), \renp_{ik(n)} = 1]\nonumber \\
& = &d(n)^2 p_1 \P[\renp_{jk(n)}=1\given N_i \geq \xi_0 k(n),
\renp_{ik(n)}=1] p_1 u_{ik(n)} \nonumber \\
& = & d(n)^2 p_1^2 u_{ik(n)} \P[\renp_{jk(n)}=1\given N_i \geq \xi_0 k(n),
\renp_{ik(n)}=1] \label{eq:summand}
\eea
Summing (\ref{eq:summand}) over $i < j$ shows that
$\sum_{i=1}^{n/k(n)} \sum_{j=i+1}^{n/k(n)} \Q_\bias[G_i^n \cap G_j^n]$ equals
\be d(n)^2p_1^2 \sum_{i=1}^{n/k(n)}u_{ik(n)} \sum_{j=i+1}^{n/k(n)}
\P[\renp_{jk(n)}=1 \given  N_i \geq \xi_0 k(n),
\renp_{ik(n)}=1]\, . \label{eq:condexps}
\ee
Let $\sigma= (i+1)k(n) - \tau_i$.  For $m=n/k(n)$, write
\be
\sum_{j=i+1}^{m} \P [ \renp_{jk(n)} = 1 \given  
N_i \geq \xi_0 k(n),\renp_{ik(n)} = 1, \sigma ] 
\label{eq:uselemma2} 
\ee
as
\be \label{eq:cetau}
\E[ \sum_{j=i+1}^{m} \renp_{jk(n)} \given 
\tau_i < (i+1)k(n) \, , \, \sigma ] 
\,.
\ee
Then observe that (\ref{eq:cetau}) is bounded above
by
\be
u_\sigma + u_{\sigma+k(n)} + \cdots + 
u_{\sigma+ m k(m)} \,.
\label{eq:uselemma}
\ee
We can apply Lemma \ref{lem:lemren} to bound
(\ref{eq:uselemma}) above by
$
\sum_{j=0}^m u_{jk(n)}
$
.
To summarize,
\be
\sum_{j=i+1}^{m} \P [ \renp_{jk(n)} = 1 \given  
N_i \geq \xi_0 k(n),\renp_{ik(n)} = 1, \sigma ] 
\; \leq \; 
\sum_{j=0}^m u_{jk(n)} \,.
\label{eq:tba}
\ee
Taking expectation over $\sigma$ in (\ref{eq:tba}),
and then plugging into
(\ref{eq:condexps}) shows that
\bea
\sum_{i=1}^{n/k(n)} \sum_{j=i+1}^{n/k(n)} \Q_\bias[G_i^n \cap G_j^n]
& \leq & 
d(n)^2p_1^2 \sum_{i=1}^{n/k(n)}u_{ik(n)}\sum_{j=0}^{n/k(n)}u_{jk(n)}
\nonumber \\
& \leq & \left( \E_\bias D_n \right)^2 + u_0 \E_\bias D_n\, ,
\label{eq:cltodone} 
\eea
where we have used the expression (\ref{eq:edn1}) for $\E_\bias D_n$.
Finally, using (\ref{eq:cltodone}) in (\ref{eq:fpsm}) yields 
that
$$
\E_\bias[D_n^2] \leq C_3 (\E_\bias[D_n])^2.$$
Now, we have, as in the proof of
Proposition \ref{prop:cntex}, that
$$\Q_\bias[ \limsup \{ D_n > 0 \} ] \geq \limsup_{n \rightarrow
  \infty} \Q_\bias[
D_n > 0 ] \geq C_3^{-1} > 0.$$
Using Lemma \ref{lem:zo} shows that
$\Q_\bias[\limsup\{D_n > 0\}] =1$.  That is, the events
$\bigcup_{i=1}^{n/k(n)}G^n_i$ happen infinitely often.  But since a
good run is also an observed run, also the events
$$\{ \exists j, \, 1\leq j \leq n/k(n) \mbox{ with } R_{jk(n)} \geq
k(n)\}$$
happen infinitely often.  But, if $R_{jk(n)} \geq k(n)$, then
certainly $R_{jk(n)} \geq k(jk(n))$.  Thus, in fact the events
$$\{ \exists j \geq k(n) \mbox{ with } R_j \geq k(j) \}$$
happen
infinitely often. That is
$$\Q_\bias[ R_n \geq k(c,n) \, i.\ o.\ ] = 1.$$
We conclude that $\c \leq
\ls(X)$.  \qed

\section{Linear estimators work when $u_n$ are
 not square-summable.} \label{sec:le}
Before stating and proving a generalization of Theorem  \ref{thm:recon},
we indicate how a weak form of that theorem may be derived
by rather soft considerations; these motivated the more concrete
arguments in our proof of Theorem \ref{thm:recongen} below.
In the setting of Theorem  \ref{thm:recon}, let
%
$${\cal T}_n \; \deq \; {\sum_{i=1}^n u_i X_i \over \sum_{i=1}^n u_i^2} \,.$$
It is not hard to verify that
$\E_\bias {\cal T}_n = \bias$ and $\sup_n \var_\bias({\cal T}_n) < \infty$.
Since $\{{\cal T}_n\}$ is a bounded sequence in $L^2(\mu_\bias)$, it
has an $L^2$-weakly convergent subsequence.  Because the limit ${\cal T}$ of
this subsequence must be a tail function, ${\cal T}=\bias$ a.s.  Finally,
standard results of functional analysis imply that
there exists a sequence of convex combinations of the estimators ${\cal T}_n$
 that tends to $\bias$ in $L^2(\mu_\bias)$ and a.s.

The disadvantage of this approach is that the convergent subsequence
and the convex combinations used
may depend on $\bias$; thus the argument sketched above
only works for fixed $\bias$.  The
proof of Theorem \ref{thm:recongen} below provides an explicit sequence
of estimators not depending on $\bias$.

We return to the general setting described in Section \ref{sec:defns}.
A collection $\Psi$ of bounded Borel functions on $\R$ is called a
{\em determining class} if $\mu = \nu$
whenever $\int_\R \psi d\mu = \int_\R
\psi d\nu$ for all  $\psi \in \Psi$.

The following theorem generalizes Theorem \ref{thm:recon}.
\begin{thm} \label{thm:recongen}
If $\sum_{k=0}^{\infty} u_k^2 = \infty$, then for any bounded Borel
function $\psi$, there
exists a sequence of functions $h_N:\R^N \to \R$ with the following
property:
\begin{quotation}
    \noindent
    for any probability measure $\eta$ on $\R$,
    we have  
    $$
    h_N(X_1,\ldots,X_{N}) \; \rightarrow \; \int \psi d\measu
    \quad \mbox{a.s. with respect to } \mu_\measu \,.
    $$
    \end{quotation}
\end{thm}
Thus the assumptions of the
theorem imply that for any countable determining class $\Psi$ of bounded Borel
functions on $\R$, a.s.\ all the integrals $\{ \int \psi d\measu \}_{\psi \in
\Psi}$ can be computed from the observations $X$, and hence a.s.\ the
measure $\eta$ can be reconstructed from the observations.

\noindent
{\sc Proof.}  Fix $\psi \in \Psi$, and assume for now that
$\measn(\psi)=\int_{\R}\psi d\alpha=0$. Without loss of generality, 
assume that $\| \psi \|_\infty \leq 1$. Define
$$w(n) = w_{n} \deq \sum_{i=0}^{n}u_{i}^{2}, \; \mbox{and} \; w(m,n) \deq
\sum_{i=m+1}^{n}u_{i}^{2}.$$
For each pair $m_i < n_i$, let
$$L_i = L_i(\psi) = \frac{1}{w(m_{i},n_{i})}\sum_{j=m_i+1}^{n_i}u_j
\psi(X_j).$$
Let $\{\ep_j\}$ be any sequence of positive numbers.  We
will inductively define $\{m_i\},\{n_i\}$ with $m_i < n_i$, so that
\be \label{eq:seqdef} w(m_{i},n_{i}) \geq w(m_{i}) \, \mbox{for all }
i, \; \mbox{and } \; \cov(L_{i},L_{j}) \leq \ep_{i} \, \mbox{for all }
j>i \, .
\ee
We now show how to define $(m_{i+1},n_{i+1})$, given $n_{i}$, so that
(\ref{eq:seqdef}) 
is satisfied.
Observe that
\bea \cov(L_i,L_\ell) &=& \frac{ \sum_{k=m_i+1}^{n_i}
  \sum_{s=m_\ell +1}^{n_\ell} u_k u_s \measu(\psi)^2 (u_ku_{s-k} -
  u_ku_s)}{w(m_{i},n_i)w(m_{\ell},n_\ell)} \nonumber \\
& = & \frac{\eta(\psi)^{2}}{w(m_{i},n_{i})w(m_{\ell},n_{\ell})}
\sum_{k=m_{i}+1}^{n_{i}} u_{k}^{2} \left(
  \sum_{s=m_{\ell}+1}^{n_{\ell}} u_{s}u_{s-k}-u_{s}^{2}\right).
  \label{eq:covbnd1}
  \eea
Fix $k$, and write $m,n$ for $m_\ell,n_\ell$ respectively.
We claim that
\be \label{eq:leclaim} \sum_{m+1}^n u_s u_{s-k}-u_s^2 \; \leq \; k \, .\ee
Assume that $\sum_{m+1}^n u_s u_{s-k}-u_s^2 > 0$; if not
(\ref{eq:leclaim}) is trivial.  Applying the inequality $a-b \leq
(a^2-b^2)/b$, valid for $b\leq a$, yields
  \be
  \sum_{s=m+1}^{n}u_{s}u_{s-k} - u_{s}^{2} \leq
  \frac{(\sum_{s=m+1}^{n}u_{s}u_{s-k})^{2} - w(m,n)^{2}}{w(m,n)}
\label{eq:precs}. \ee
Then applying Cauchy-Schwarz to the right-hand side of
(\ref{eq:precs}) bounds it by 
\bea
\frac{w(m,n)w(m-k,n-k) - w(m,n)^{2}}{w(m,n)} 
& \leq & w(m-k,n)-w(m,n) \nonumber \\
& = & w(m-k,m) \nonumber \\
& \leq & k \,, \nonumber
\eea
establishing (\ref{eq:leclaim}).
Using the bound (\ref{eq:leclaim})
in (\ref{eq:covbnd1}), and recalling that $|\psi| \leq
1$, yields \be \cov(L_i,L_\ell) \leq
\frac{1}{w(m_i,n_i)w(m_{\ell},n_{\ell})}\sum_{k=m_i+1}^{n_i}
u_k^2k \leq 
\frac{n_{i}}{w(m_{\ell},n_{\ell})} \; . \label{eq:covbnd3} \ee Pick
$m_{i+1}$ large enough so that \be \label{eq:midef} w(m_{i+1}) \geq
\frac{n_{i}}{\ep_{i}} \; ,\ee and let $n_{i+1} \deq \inf\{t:
w(m_{i+1},t) \geq w(m_{i+1}) \}$.  Then for any $\ell \geq i+1$, since
$w(m_{\ell},n_{\ell}) \geq w(m_{\ell}) \geq w(m_{i+1})$,
(\ref{eq:midef}) and (\ref{eq:covbnd3}) yield that $\cov(L_i,L_\ell)
\leq \ep_i.$

Observe that $\E[ L_i ] \; = \; \measu(\psi)$, and 
\bea 
\E\left[ \left(\sum_{j=m_i+1}^{n_i}u_j\psi(X_j)\right)^2 \right] 
& = & 2 \measu(\psi)^2 \sum_{j=m_i+1}^{n_i}\sum_{k=j+1}^{n_i}u_j u_k u_j
u_{k-j} \nonumber \\
&& + \; \sum_{j=m_i+1}^{n_i}\E[\psi(X_j)^2 ] \, u_j^2 \nonumber \\
& \leq &
\|\psi\|^2_\infty\left\{ 2\sum_{j=m_{i}+1}^{n_{i}}u_{j}^{2}
\sum_{k=j+1}^{n_{i}}u_{k}u_{k-j}
+ w(m_{i},n_{i}) \right\} \,.
\label{eq:sm1}
\eea Fix $i$, let $m=m_i,n=n_i$.  For $j$ fixed, using Cauchy-Schwarz
yields \be \sum_{k=j+1}^n u_k u_{k-j} \; \leq \; \sqrt{
  w(j,n)w_{n-j}} \; \leq \; w_{n} \, . \label{eq:cs1} \ee Plugging
(\ref{eq:cs1}) into (\ref{eq:sm1}), and recalling that 
$\|\psi\|_\infty < 1$,
gives that \be \E\left[
  \left(\sum_{j=m_i+1}^{n_i}u_j\psi(X_j)\right)^2 \right] \; \leq \;
2 w_{n_i}^2 + w_{n_i}\; . \nonumber \ee
Thus,
$$
\E[ L_i^2 ] \; \leq \; \frac{2 w_{n_i}^2 +
  w_{n_i}}{w_{n_i}^2 / 4} \; = \; 8 + \frac{4}{w_{n_i}}
\; \leq \; B \; .
$$
Choosing, for example, $\ep_{i} = i^{-3}$, one can apply the strong
law for weakly correlated random variables (see Theorem A
in section 37 of \cite{Lo}), to get that
\be \label{eq:converge} G_n(\psi) \deq \frac{1}{n}\sum_{i=1}^n
L_i(\psi) \rightarrow \measu(\psi) \; a.\ s.  \ee
For general $\psi$, define $H_n(\psi) = G_n(\psi - \measn(\psi)) +
\measn(\psi)$.   From (\ref{eq:converge}), it follows that
\be H_n(\psi) \rightarrow \eta(\psi - \measn(\psi)) + \measn(\psi) =
\eta(\psi) \, .\nonumber \ee
To finish the proof, define $h_N(X_1,\ldots,X_N) \deq H_{k(N)}(\psi)$,
where $k(N)$ is the largest integer $k$ such that $n_k \leq N$.
\qed

%
\section{Quenched Large Deviations Criterion.} \label{sec:qmgf}
Recall that $\rho_n = {d\mu_{\measu}\over  d\mu_{\measn}}_{\vline \G_n}$, the
density of the measure $\mu_\measu$ restricted to $\G_n$ with respect
to the measure $\mu_\measn$ restricted to $\G_n$.

We make the additional assumption that
\be r \; = \; \int_\R \left({d\measu \over d\measn}\right)^2 d\measn
\; = \; \int_\R {d\measu \over d\measn} d \measu \; < \; \infty.
\label{eq:rdefn}
\ee
For two binary sequences $\rent,\rent'$,
define 
$J(\rent,\rent') = | \{n: \rent_n = \rent'_n
= 1\}|$, the number of joint renewals.
\begin{lem} \label{lem:qmgf}
If $\E [r^{J(\renp,\renp')}\given \renp ] < \infty$, then
$\mu_\measu \ll
\mu_\measn$.
\end{lem}
{\sc Proof.}  Let $x(\nt,\st,\rent)_n = \st_n\rent_n +
\nt_n(1-\rent_n)$. We have
\be
\E_{\Q_{\eta}} [ \rho_n(X) \given \renp=\rent]
\; = \; \int\limits_{\R^{\infty}}\int\limits_{\R^{\infty}}
\rho_n (x(\nt,\st,\rent)) d\measn^{\infty}(\nt)
d\measu^{\infty}(\st) \,,
\label{eq:qre}
\ee
and expanding $\rho_n$ shows that (\ref{eq:qre}) equals
\be
\int\limits_{\R^{\infty}}\int\limits_{\R^{\infty}}
\int\limits_\Upsilon \prod_{i=1}^n \left[
{d\measu \over d\measn}(x(\nt,\st,\rent)_i)\rent_i' +
1-\rent_i' \right] d\P(\rent')d\measn^{\infty}(\nt)d\measu^{\infty}(\st)
\,. \label{eq:lem1} 
\ee
Using Fubini's Theorem and the independence of coordinates under
product measure, (\ref{eq:lem1}) is equal to
\be \label{eq:laterplug}
\int_\Upsilon \prod_{i=1}^n \int_\R \int_\R
\left[ {d\measu \over d\measn}(x(\nt,\st,\rent)_i)\rent_i' +
1-\rent_i' \right]d\measn(\nt)d\measu(\st)d\P(\rent').
\ee
If
$$
I \; \deq \; \int_\R \int_\R
\left[ {d\measu \over d\measn}(x(\nt,\st,\rent)_i)\rent_i' +
1-\rent_i' \right]d\measn(\nt)d\measu(\st)
\,,
$$
then we have that
\be \label{eq:cases}
I \; = \; \left\{ \begin{array}{ll}
1 & \mbox{if }\rent_i'=0 \\
\int {d\measu\over d\measn}(\nt) d\measn(\nt)=1
& \mbox{if } \rent_i'=1,\rent_i =0\\
\int {d\measu\over d\measn}(\st) d\measu(\st)=r
& \mbox{if } \rent'=1, \rent_i=1
\end{array} \right. 
\ee
Plugging (\ref{eq:cases}) into (\ref{eq:laterplug}), we get that
\bean
\E_{\Q_\measu} [ \rho_n(X) \given \renp=\rent] & = &
 \int_\Upsilon \prod_{i=1}^n r^{\rent_i\rent_i'}d\P(\rent') \\
& \leq & 
\int_\Upsilon r^{J(\rent,\rent')} d\P(\rent') \\
& = & \E[ r^{J(\renp,\renp')} \given \renp=\rent] \,.
\eean
Applying Fatou's Lemma, we infer that $\E{\Q_\measu}[\rho(X) \given \renp] <
 \infty$, whence
$$\Q_\measu[ \rho(X) < \infty] \; = \; 1 \,.$$
The Lebesgue Decomposition
(\ref{eq:lebdec}) implies that $\mu_\measu \ll \mu_\measn$.  \qed

\section{Absence of Phase Transition in Almost Transient Case}
\label{sec:example}
In this section, we apply the quenched moment generating function
criterion established in the previous section.

Let $N[m,n]$ be the
number of renewals in the interval $[m,n]$, and write $N_m = N[0,m]$.
Let $U_m = U(m) = \E
N_m = \sum_{k=0}^m u_k$.
\begin{lem} For any integer $A \ge 1$, we have
 $\P[ N_m \geq A e U_m] \leq e^{-A}$ .\end{lem}
{\sc Proof.} For $A=1$, the inequality follows from Markov's
inequality.  Assume it holds for $A-1$.  On the event $E$ that $N_m
\geq (A-1) e U_m$, define $\tau$ as the time of the $\lceil (A-1) e
U_m \rceil^{th}$ renewal.  Then
$$
\P[ N_m \geq A e U_m \given E ] \leq \P [ N[\tau, m] \geq e U_m|E ]
\leq \P[ N_m \geq e U_m] \leq e^{-1}.
$$
Consequently,
$$ \P[N_m \geq A e U_m ] \leq \P[N_m \geq A e U_m | E]e^{-(A-1)}
\leq e^{-A}.$$ \qed
\begin{thm} \label{thm:nophase}
Suppose that the renewal probabilities $\{u_n\}$ satisfy 
$$ U(e^{k})  \; = \; o(k/\log k) \,,$$ 
and also $u_k \leq C_2k^{-1}$.
If $\measu \ll \measn$ and ${d\measu\over d\measn} \in L^2(\measn)$,
then $\mu_\measu \ll \mu_\measn.$
\end{thm}
{\sc Proof.}
In this proof the probability space will always be $\ren^2$, endowed
with the product measure $\P^2$, where $\P$ is the renewal
probability measure.
Let 
$$
J[m,n] \; = \; |\{n \leq k \leq m: \renp_k = \renp'_k = 1\}|
$$ 
be the
number of {\em joint renewals} in the interval $[m,n]$.

First we show that
\be \forall \; C, \quad T_1 + \cdots + T_k \geq e^{Ck} \; \;
\mbox{ eventually}. \ee
Observe that
\be
\P[T_1 + \cdots + T_k \leq e^{Ck} ] = \P[ N(e^{Ck}) \geq k ]
\leq  \exp\left( - \frac{k}{eU(e^{Ck})} \right) \, .
\label{eq:bndn}
\ee
Our assumption guarantees that
$k/eU(e^{Ck}) \ge 2\log k$ eventually, and
hence the right-hand side of (\ref{eq:bndn}) is summable.  Consequently, for
almost all $\renp$, there is an integer $M=M(\renp)$
such that  $\sum_{j=1}^k T_j > e^{Ck}$ for all $k > M(\renp)$.
Equivalently, $N[0,\exp(Ck)] < k$ when $k > M$.
To use Lemma \ref{lem:qmgf}, it suffices to show that $$\sum_n s^n
\P[J[0,n] \geq n \given \renp ] < \infty  \; \mbox{ a.s., for all real
  } s.$$
We have
\be
\sum_n s^n \P[J[0,n] \geq n \given \renp ]  \; \leq \;
C_2(\renp) + \sum_{n=M}^{\infty} s^n \P[J(e^{Cn},\infty) \geq 1 \given \renp ]
\label{eq:quen} \ee
Observe that
\be \E[ J(e^{Cn},\infty)] =
\sum_{k=\exp(Cn)}^{\infty} u_k^2 \leq C_3 e^{-Cn}, \ee
since we have assumed that $u_n \leq C_2n^{-1}$.
Thus the expectation of the sum on the right in (\ref{eq:quen}),
for $C$ large enough, is finite.  Thus the sum is finite
$\renp$-almost surely, so the conditions of Lemma \ref{lem:qmgf} are
satisfied.  We conclude that $\mu_\measu \ll \mu_\measn$. \qed

We now discuss examples of Markov chains which satisfy the
hypothesis of Theorem \ref{thm:nophase}.

\begin{lem} \label{lem:comp} Given two Markov chains with transition
matrices
$P,P'$ on state spaces ${\cal X}$ and ${\cal Y}$ with distinguished
states $x_0, y_0$ respectively,  construct a new
chain $\mc = (X,Y)$ on ${\cal X} \times {\cal Y}$ with transition matrix
$$Q((x_1,y_1),(x_2,y_2)) = \left\{ \begin{array}{ll}
 P(x_1,x_2)P'(y_1,y_2) & \mbox{if } y_1 = 0 \\
 P'(y_1,y_2) &\mbox{if } y_1 \neq 0, \; x_1 = x_2
 \end{array} \right. \, .$$
Let $A(s) = \sum_{n=1}^{\infty} f_n s^n$ be the moment generating
function for the distribution of the time of first return to $x_0$ for the
chain with transitions $P$, and let $B(s)$ be the corresponding generating
function but
for the chain $P'$ and state $y_0$.  Then the generating function for the
distribution
of the time of the first return of $\mc$ to $(x_0,y_0)$ is the
composition $A \circ B$.
\end{lem}
{\sc Proof.}
Let $S_1, S_2, \ldots$ be the times of successive visits of $\mc$ to ${\cal X}
\times \{y_o\}$, and $T_k = S_k - S_{k-1}$.  Observe that $Y$ is a Markov
chain with transition matrix $P'$, so $\{T_k\}$ has the
distribution of return times to $y_0$ for the chain $P'$.

Let $\tau = \inf\{n\geq 1: X_{S_n} = x_0\}$.
Note that $\{X_{S_n}\}_{n=0}^{\infty}$ is a Markov chain with
transition matrix $P$, independent of $\{T_n\}$.  Hence $\tau$ is
independent of $\{T_n\}$, and
$$T = T_1 + \cdots + T_{\tau}$$ is the time of the first return of
$\mc$ to $(x_0,y_0)$.  A standard calculation (see, for example,
XII.1 in \cite{F1}) yields that the generating function $\E s^T$ is $A
\circ B$. \qed

Let $F,U$ be the moment generating functions for the sequences $\{f_n\}$
and $\{u_n\}$ respectively.
Define $L:(0,\infty) \rightarrow (1,\infty)$ by $L(y) = 1-
{1\over y}$, and note that $F=L\circ U$. Denote
$W(y) = U\circ L(y) = L^{-1}\circ F\circ L$.
When $F_3 = F_1 \circ F_2$, it follows that $W_3 = W_1\circ W_2$.

We  use the following Tauberian theorem from \cite[Theorem 2.4.3]{La}:
\begin{prop} \label{prop:taub}
Let $\{a_n\}$ be a sequence of non-negative reals,
$A(s) = \sum_{n=0}^\infty a_ns^n$ its
generating function, $W(y) \deq A(1 - y^{-1})$, $\alpha \geq 0$ a constant,
and $\ell$ a slowly varying function.  The following
are equivalent:
\begin{description}
\item{(i)} $A(s) \asymp (1-s)^{-\alpha}\ell((1-s)^{-1})$ for $s<1$
near 1.
\item{(ii)} $W(y) \asymp y^\alpha \ell(y)$ for large $y$.
\item{(iii)} $A_n = \sum_{k=0}^n a_k \asymp n^\alpha \ell(n)$ \, .
\end{description}
\end{prop}

We now exhibit Markov chains with no phase transition.
\begin{prop} \label{prop:nophaseex} There is a Markov chain that
satisfies
$U_n \asymp \log \log n$,
  and $u_n \leq Cn^{-1}$.
\end{prop}
{\sc Proof.}
For simple random walk on $\Z^2$, we have 
$$
U(s) \asymp
\sum_{n=1}^{\infty}n^{-1}s^{2n} = -\log(1-s^2)
\,.
$$ 
Thus, $W(y) \asymp \log
y$.  Consequently, $W \circ W (y) \asymp \log \log (y)$ corresponds
to the chain in Lemma \ref{lem:comp} with both $P$ and $P'$ the
transition matrices for simple random walk on $\Z^2$.
Proposition \ref{prop:taub} implies that $U_n \asymp \log\log n$.
Finally, 
$$u_n \leq \P[X_n = 0] \leq Cn^{-1} \,,
$$
since $X$ is a simple random walk on $\Z^2$. \qed

In conjunction with Theorem \ref{thm:nophase},
this establishes Theorem \ref{thm:nopt}.

Lemma \ref{lem:comp} can be applied to construct Markov chains obeying
the hypotheses of Proposition \ref{prop:cntex} and Theorem
\ref{thm:reconstr}.  Take as the chains $X$ and $Y$ the simple random
walk on $\Z$.  The moment generating function $U_{[1]}$ for the return
probabilities $u_n$ of the simple random
walk is given by $U_{[1]}(s) =
(1-s^2)^{-1/2}$ (see XIII.4 in \cite{F1}).  Then
$W_{[1]}(y) = U_{[1]}\circ L(y)= ({y
\over 2- y^{-1}})^{1/2}$ satisfies $W_{[1]}(y) \sim (y/2)^{1/2}$ as $y
\rightarrow \infty$.  Hence $W(y) = W_{[1]} \circ W_ {[1]}(y) \asymp y^{1/4}$,
and by Proposition \ref{prop:taub}, $U_n \asymp n^{1/4}$.

The last example is closely related to the work of
Gerl in \cite{G}. He considered certain
``lexicographic spanning trees'' ${\cal T}_d$ in  $\Z^d$, where the
path from the origin to a lattice point $(x_1,\ldots,x_d)$
consists of at most $d$ straight line segments, going through the points
$(x_1,\ldots,x_k,0,\ldots,0)$ for $k=1,\ldots,d$ in order. 
Gerl  showed that for $d \ge 2$, the return probabilities of
simple random walk on ${\cal T}_d$ satisfy
 $u_{2n} \asymp n^{2^{-d}-1}$; after
introducing delays, this provides further
examples of Markov chains with a phase transition ($0<\theta_c<1$).

\section{Absence of Phase Transition in $\Z^2$.}
\label{sec:z2}
%
%
The results in \cite{HK} (as summarized in
Theorem B of Section \ref{sec:intro}) show that for simple random walk on 
$\Z^2$, which moves in each step to a uniformly chosen neighbor, 
the measures $\mu_\bias$ and $\mu_0$ are mutually
absolutely continuous for all $\bias$.
The argument does not extend 
to Markov chains which are small perturbations of this walk.  For 
example, if the walk is allowed to remain at its current position
with some probability, the asymptotic behavior of $\{u_n\}$ is not
altered, but Theorem B does not resolve whether $\mu_\bias \ll 
\mu_0$ always.  In this section, we show that for any Markov chain 
with return probabilities that satisfy
$u_n =O( n^{-1})$,
the measures $\mu_\bias$ and $\mu_0$ are mutually absolutely 
continuous.

Recall that $T$ is the time of the first renewal, and $T_1,T_2,\ldots$
are \iid copies of $T$. Also, $S_n = \sum_{j=1}^n T_j$ denotes the
time of the $n^{th}$ renewal.  Recall from before 
that $\renp_n$ is the indicator
of a renewal at time $n$, hence
$$
\{ S_n = k  \mbox{ for some } n \geq 1\} \; = \; \{ \renp_k = 1 \}
\,.
$$
Let $S_n'$ and $T_n'$ denote the renewal times and inter-renewal
times of another independent renewal process.  
Recall that $J$ is the total number of simultaneous renewals:
$J = \sum_{k=0}^\infty \renp_k \renp'_k$.
If $\rF_k$ is the
sigma-field generated by $\{ T_j : 1 \leq j \leq k \}$, then
define
\be \label{eq:qn}
q_n \; = \; \P[ J \geq n \given \renp]
\; = \; \P \Big[\; | \, \{ (i,j) : S_i = S_j' \} \, | \; \geq n 
 \Big|  \rF_\infty \Big] \,.
\ee
In this section, we prove the following:
\begin{thm} \label{thm:Z2}
When $u_n = O( n^{-1})$, the sequence $\{ q_n \}$ defined in
(\ref{eq:qn}) decays 
faster than exponentially almost
surely, that is, 
$$
n^{-1} \log q_n \to -\infty \mbox{ almost surely.}
$$
Consequently, the quenched large deviations criterion Lemma
\ref{lem:qmgf} implies that
if $\measu \ll \measn$ and ${d\measu\over d\measn} \in L^2(\measn)$,
then $\mu_\measu \ll \mu_\measn.$
\end{thm}

We start by observing that the assumption
$u_n \leq c_1 / n$ implies a bound for tails of
the inter-renewal times:  
\begin{equation} \label{eq:logtail}
\exists  c_2>0 \quad \P [\log T \geq t] \geq c_2 t^{-1} .
\end{equation}
Indeed, by considering the last renewal before time $(1+a)n$,
\begin{eqnarray*}
1&=&\sum_{k=0}^{(1+a)n-1} u_k \P[T \ge (1+a)n-k] \\[1ex]
& \le & \sum_{k=0}^{an} u_k \P[T \ge n] + \sum_{k=an+1}^{(1+a)n} u_k   \\[1ex]
& \le  & (2+c_1 \log an)  \P[T \ge n] +2c_1 \log \frac{1+a}{a} \,.
\end{eqnarray*}
Choosing $a$ large yields (\ref{eq:logtail}).

Let  $\omega (n)$ be any function going to
infinity, and denote
$$
 m(n) := n \log n \omega^2 (n) \,.
$$
Below, we will often
write simply $m$ for $m(n)$.

{}From  (\ref{eq:logtail})
it follows that 
$$\P [S_{m(n)} \leq e^{n \omega (n)} ] \leq
   \Big(1 - \frac{c} {n \omega (n)}\Big)^{m(n)} \leq
   n^{-c \omega (n)}.$$
This is summable, so by Borel-Cantelli, 
\begin{equation} \label{eq:far}
n^{-1} \log S_{m(n)} \to \infty
\end{equation}
almost surely.

\proof{Proof of Theorem~\ref{thm:Z2}} 
Define 
the random variables
$$
J_m \; \deq \;
| \{ (i,j) \; : \; i > m, j \geq 1 \mbox{ and } S_i = S_j' \; \} |
$$
and let
$
Q_m \deq 
\P[ J_{m} \geq 1 \given \rF_\infty ]
$
.

Let
$$
r_n \; \deq \; \P\Big[ \; | \, \{ (i,j) \; : \;
i \leq m(n) \mbox{ and } S_i = S'_j \} \,| \; \geq n \; \Big| \; 
\rF_\infty \Big] \,.
$$ 
Clearly,
\be \label{eq:qdecomp}
q_n 
\leq  Q_{m(n)} + r_n \,.
\ee
Write
$
Q^*_m \; \deq \; \E [ Q_m \given \rF_m ] \; = \; 
\P[ J_{m} \geq 1
\given \rF_m] 
$
.  
Then
\begin{eqnarray*}
Q^*_{m(n)} \; \leq \; \E[ J_{m(n)} \given \rF_{m(n)} ] 
& \leq  & \; \sum_{k=1}^\infty u_k u_{k + S_{m(n)}} \\
& \leq & \; \sum_{k=1}^\infty
   {c_1 \over k} {c_1 \over k + S_{m(n)}} 
\; \leq \; c_3 {\log S_{m(n)} \over S_{m(n)}} \,.
\end{eqnarray*}
By~(\ref{eq:far}), we see that $n^{-1} \log Q_{m(n)}^* \to
-\infty$ almost surely. 

Since $Q_{m(n)}^* = \E [Q_{m(n)} | \rF_{m(n)}]$, we see that
$\P [Q_{m(n)} \geq 2^n Q_{m(n)}^*] \leq 2^{-n}$, hence 
$$Q_{m(n)} \geq 2^n Q_{m(n)}^* \mbox{ finitely often}$$
and it follows that 
$
n^{-1} \log Q_{m(n)} \to -\infty 
$.
It therefore suffices by (\ref{eq:qdecomp}) 
to show that 
\be \label{eq:subexp}
\frac{\log r_n}{n} \to -\infty
\mbox{ almost surely.}
\ee
Let $[m(n)] \deq \{ 1 , 2 , \ldots , m(n) \}$.
We can bound $r_n$ above by
\be \label{eq:rnsum}
\sum_{\stackrel{A\subset [m(n)]}{|A|=n}}
\P[ \forall i \in A, \exists j \geq 1 \mbox{ so that }
S'_j = S_i \given  \rF_{\infty}\, ] 
\; \leq \; {m(n) \choose n} R_n \,, 
\ee
where
$$
R_n \; \deq \;
\max_{\stackrel{A\subset [m(n)]}{|A|=n}}
\P[ \forall i \in A, \exists j \geq 1 \mbox{ so that }
S'_j = S_i \; \given \; \rF_{\infty}\, ] \,.
$$
We can conclude that
\be \label{eq:logineq}
\log r_n \; \leq \; \log {m(n) \choose n} + \log R_n 
\,.
\ee
Notice that ${m(n) \choose n} = e^{O(n \log \log n)}$ when $\omega (n)$ is no
more than ${\rm polylog}\; n$;  for convenience, we assume
throughout that $\omega^2 (n) = o(\log n)$. 
Hence, 
if we can show that 
\begin{equation} \label{eq:R}
{\log R_n \over n \log \log n} \to -\infty
\mbox{ almost surely,}
\end{equation}
then by (\ref{eq:logineq}), it must be that (\ref{eq:subexp}) holds.

For any $n$-element set $A \subset [m(n)]$, we use the following 
notation:
\begin{itemize}
    \item $A \; = \; \{x_1 < x_2 < \cdots < x_n \}$, and $m' \deq x_n$.
    \item For any $k \leq m'$, let 
    $I(k)$ be the set of indices $i$ such that $\{ T_i \}_{i \in I(k)}$
    are the $k$ largest inter-renewal times among 
    $\{ T_i \}_{i \leq m'}$. 
    \item For $i \leq n$, let 
    ${M(A,i)} \deq \max\{ T_j : x_{i-1}+1 \leq j \leq x_i \}$.
\end{itemize}
We have
$$
\P[ \forall x_i \in A, \exists j \geq 1 \mbox{ so that }
S'_j = S_{x_i} \; \given \; \rF_{\infty}\, ] \; = \;
\prod_{i=1}^n
u_{S_{x_i} - S_{x_{i-1}}}\,,$$
where $x_0=0$ and $S_0 \deq 0$.  Recalling that 
$u_n \leq c_1 / n$, we may bound the \rhs above by
\begin{equation} \label{eq:R(A)}
\prod_{i=1}^n \frac{c_1}{S_{x_i}-S_{x_{i-1}}} \; = \;
\prod_{i=1}^n \frac{c_1}{\sum_{j=x_{i-1}+1}^{x_i}T_j} 
\; \leq \;  R(A) \deq \prod_{i=1}^n {c_1 \over {M(A,i)}} \,.
\ee
To summarize, we have
\be \label{eq:sumrn}
R_n \; \leq \; \max_{\stackrel{A\subset [m(n)]}{|A|=n}}
R(A) \; = \; \max_{\stackrel{A\subset [m(n)]}{|A|=n}}
\prod_{i=1}^n {c_1 \over {M(A,i)}}
\,.
\ee
To see where this is going, compute what happens when $A = [n]$.
{}From the tail behavior of
$T$, we know that 
$$
\liminf_{n \rightarrow \infty}{\log R([n]) \over n \log n} > 0 \,.
$$
To
establish~(\ref{eq:R}), we need something like this for $R_n$
instead of $R([n])$.

In what follows, 
$k_0(n) \deq 10 (\log n \omega(n) )^2$.
\begin{lem} \label{lem:big}
Almost surely, there is some (random) $N$ so that if
$n > N$, then 
    for all 
    $n$-element sets $A \subseteq [m(n)]$, providing  
    $k$ satisfies $m' \geq k > k_0(n)$, at least
    $k n / (6 m' \log \log n)$ 
    values of $i$ satisfy 
    $M(A,i) \in \{T_j \; : \;j \in I(k) \}$.
\end{lem}

Assuming this for the moment, we finish the proof of the theorem.
The following summation by parts principle will be needed.
\begin{lem} \label{lem:parts}
Let $\H(k)$ be the $k$ largest values in a given finite set $\H$ of
positive real numbers. Suppose another set $\H'$ contains at least
$\ep k$ members of $\H(k)$ for every $k_0 < k \leq |\H|$.  Then 
$$\sum_{\h \in H'} \h \; \geq \; \ep 
\sum_{\h \in \H \setminus \H(k_0) } \h
\,.
$$
\end{lem}
\proof{Proof of Lemma \ref{lem:parts}}
Let $\H=\{\h_j, j=1,...,N\}$
in decreasing order and let $\h_{N+1}=0$ for convenience.
Write 
$$
f(j) \deq \one_{\{ \h_j \in \H' \}} \,,
\mbox{ and let }
F(k)=f(1)+...+f(k)\,. 
$$
Then
\begin{eqnarray*}
\sum_{j=1}^N f(j)\h_j & = & 
\sum_{j=1}^N (F(j)-F(j-1)) h_j 
\; = \; 
\sum_{k=1}^N F(k) (\h_k-\h_{k+1}) \\
& \geq & \sum_{k = k_0 + 1}^N F(k) (\h_k - \h_{k+1}) 
\; \geq \; \sum_{k = k_0 + 1}^N \epsilon k (\h_k-\h_{k+1})\\
& = & \ep\left\{ (k_0 + 1)\h_{k_0 +1} + \sum_{k= k_0 + 2}^N \h_k
\right\} 
\; \geq \; \ep \sum_{k= k_0 + 1}^N \h_k \\
& = & \ep \sum_{\h \in \H\setminus \H(k_0)} \h
\end{eqnarray*}
This proves the lemma. \qed
\endproof
\begin{lem}\label{lem:omitsum}
    Write $\{T_i\}_{i=1}^n$ in decreasing order:
    $$T_{(1)} \geq T_{(2)} \geq \cdots \geq T_{(m)}\,.$$
Then
$$
\liminf_{n \rightarrow \infty}
\frac{1}{n\log n}\sum_{i = k_0(n)+1}^n \log T_{(i)} 
> 0
\,.
$$
\end{lem}
\proof{Proof} 
It suffices to prove this lemma in the case where $u_n \asymp n^{-1}$,
because in the case where $u_n \leq c n^{-1}$, the random variables
$T_i$ stochastically dominate those in the first case.

Let $Y_i \deq \log T_i$; then $Y_i$ are \iid random variables
with tails obeying
$$ \P[Y_i \geq t ] \asymp t^{-1}
\,.$$
Write $Y_{(i)}$ for the $i^{th}$ largest among $\{Y_i\}_{i=1}^n$.
{}From \cite{CS}, it can be seen that
\be \lab{eq:sumlog}
\lim_{n \rightarrow \infty} \frac{1}{n\log n}
\left( \sum_{i=2}^{k_0(n)} Y_{(i)} - n\log \log n \right)
= 0
\,.
\ee
{}From Theorem 1 of \cite{Mo}, we can deduce that
\be \label{eq:sumwomax}
\liminf_{n \rightarrow \infty} \frac{1}{n\log n}
\sum_{i=2}^n Y_{(i)}
> 0 \,.
\ee
%

Combining (\ref{eq:sumlog}) and (\ref{eq:sumwomax}) yields
$$
\liminf_{n \rightarrow \infty} \frac{1}{n\log n}
\sum_{i=k_0(n)+1}^n Y_{(i)}
> 0 \,.
$$
\qed
\endproof

Recall that $$R(A) = \prod_{i=1}^n c M(A,i)^{-1}\,.$$
{}From 
Lemma~\ref{lem:big} we see that almost surely there exists
an $N$ so that, for all $n>N$ and $k_0(n) < k \leq m'$ 
the set $\{M(A,i) \; : \; 1 \leq i \leq n\}$
includes at least
$k n / (6 m' \log \log n)$ of the $k$ greatest values of
$\{T_j\}_{j=1}^{m'}$ 
Therefore by Lemma~\ref{lem:parts}
(applied to the logs of the denominators), we see that 
for $n>N$ and all $A \subset [m(n)]$,
$$- \log R(A) \geq {(n/m') \sum_{i=k_0(n)+1}^{m'} \log (T_{(i)} / c)
   \over (6 \log \log n)}
   \,.$$
Since $(m' \log m')^{-1} \sum_{i=k_0(n)+1}^{m'}
\log (T_i / c)$ has a nonzero liminf by Lemma \ref{lem:omitsum},
we see that
$\log \log n {\log R_n \over n \log n}$ is not going to 
zero, from which follow~(\ref{eq:R}) and the theorem. \qed
\endproof

It remains to prove Lemma~\ref{lem:big}. 
Define the event
$G_{n,m'}$ to be the event
\begin{quote}
for all $n$-element sets $A \subset [m(n)]$
with maximal element $m'$,
and $k$ obeying $m' \geq k > k_0 (n)$, at least 
$k n / (6 m' \log \log n)$ 
values of $i$ satisfy 
$$M(A,i) \in \{T_i \; : \;i \in I(k) \}\,.$$
\end{quote}
Then define $G_{n} \deq \cap_{m'=n}^{m(n)} G_{n,m'}$.  The conclusion
of Lemma~\ref{lem:big} is that
\be \label{eq:biglemcon}
\P[ G_n \mbox{ eventually} ] \; = \; 1 \,.
\ee
If we can show that
\be \label{eq:showin}
\P[ G_{n,m'}^c ] \; \leq \; n^{-3} \,,
\ee
then by summing over $m' \in [n,m(n)]$, we can conclude that
$\P[G_n^c] \leq \frac{\log n \omega^2(n)}{ n^{2}}$, and hence by Borel-Cantelli, that
(\ref{eq:biglemcon}) holds.

We prove (\ref{eq:showin}) for $m' = m$, the
argument for other values of $m$ being identical.  The values
$T_1 , T_2 , \ldots$ are exchangeable, so the set $I(k)$ is a
uniform random $k$-element subset of $[m]$ and we may restate
(\ref{eq:showin}) (with $m'=m$):

Let 
$$
I(k) \; = \; \{ r_1 < r_2 < \cdots < r_k \}
$$
be a uniform $k$-subset of $[m(n)]$; then the
event $G_{n,m}$ has the same probability as the event
$\widetilde{G}_{n,m}$, defined as
\begin{quote}
for all $n$-element sets 
$A = \{ x_1 < \cdots < x_n = m \} \subseteq [m]$
and $k$ satisfying $m \geq k > k_0(n)$, 
at least
$k n / (6 m \log \log n)$ of the intervals $[x_{i-1}+1 , x_i]$
contain an element of $I(k)$.  
\end{quote}
Equivalently, $\widetilde{G}_{n,m}$
is the event that
\begin{quote}
for all $n$-element sets 
$A = \{ x_1 < \cdots < x_n = m \} \subseteq [m]$
and $k$ satisfying $k > k_0(n)$,   
at least 
$k n / (6 m \log \log n)$ of the intervals 
$[r_i , r_{i+1} - 1] , 1 \leq i \leq k$
contains an element of $A$.
\end{quote}
Finally,
$\widetilde{G}_{n,m}$ can be rewritten again as the event
\begin{quote}
for $k$ obeying $m \geq k > k_0(n)$,
no $k n / (6 m \log \log n) - 1$ 
of the intervals $[r_i , r_{i+1} - 1]$ together contain
$n$ points.
\end{quote}
Proving the inequality (\ref{eq:showin}) is then the same as proving
that
\be \label{eq:showin2}
\P[ \widetilde{G}_{n,m} ] \; \geq \; 1 - n^{-3}
\,.
\ee
For $0 \leq j \leq k$ let $D_j$ denote $r_{j+1} - r_j$ where
$r_0 := 0$ and $r_{k+1} \deq m+1$.  For any $B \subseteq [k]$,
let $W(B)$ denote the sum $\sum_{j \in B} D_j$.  Then define
the events
$\widetilde{G}_{n,m,k}$
to be
\begin{quote}
For all sets $B \subset [k]$ with 
$|B| < k n / (m \log \log n)$,
we have
$W(B) < n$.
\end{quote}
We have that $\widetilde{G}_{n,m} = \cap_{k=k_0(n)+1}^m
\widetilde{G}_{n,m,k}$.

Set $\ep = n / m = (\log n \omega^2
(n))^{-1}$, and set $\delta = \ep / (6 \log \log n)$, so that
$$\delta \log {1 \over \delta} = {\ep \over 6} {\log (1/\ep) \over
   \log \log n} \leq {2 \ep \over 5}$$
for sufficiently large $n$.  We 
now need to use the following lemma:
\begin{lem} \label{lem:1}
Let $p(k , m , \ep , \delta)$ denote the probability that there is
some set $B$ of cardinality at most $\delta k$ such that
$W(B) \geq \ep m$.  Then for $\ep$ sufficiently small and
$\delta \log (1/\delta) \leq \ep / 5$,
$$p (k , m , \ep , \delta) \leq e^{- k \ep / 2} .$$
\end{lem}
The proof of this will be provided later.

Now applying Lemma \ref{lem:1}, 
we have that for fixed $k$ so that $m \geq k > k_0(n)$,
$$
\P[\widetilde{G}_{n,m,k}] \; \geq  \; 1 - n^{-5} \,,
$$
since $\frac{k\ep}{2} \geq n^{-5}$.
Summing over $k$ gives that
$$
\P[\widetilde{G}_{n,m}]\; \geq \;1 - n^{-3} \,.
$$

To prove Lemma \ref{lem:1}, two more lemmas are required.
\begin{lem} \label{lem:char}
Let $B \subseteq [k]$ and $W := \sum_{j \in B}
D_j$.  Then for $0 < \lambda < 1$,
\begin{equation} \label{eq:char}
\E e^{\lambda k W / m} \leq \left ( {1 \over 1 - \lambda}
   \right )^{|B|} .
\end{equation}
\end{lem}

\noindent
\proof{Proof} The collection $\{ D_j : 0 \leq j \leq k \}$
is exchangeable and is stochastically increasing in $m$.  It
follows that the conditional joint distribution of any subset of
these given the others is stochastically decreasing in the values
conditioned on, and hence that for any $B \subseteq [k]$, and
$\lambda > 0$,
\begin{equation} \label{eq:NA}
\E \exp \left ( \sum_{j \in B} D_j \right ) \leq
   \prod_{j \in B} \E \exp (D_j) = \left ( \E \exp (D_0)
   \right )^{|B|} .
\end{equation}
The distribution of $D_0$ is explicitly described by
$$
\P (D_0 \geq j) \; = \; (1 - {j \over m}) \cdots (1 - {j \over
m-k+1}) \,.
$$  Thus
$$\P (D_0 \geq j) \leq \left ( 1 - {j \over m} \right )^k
   \leq e^{-kj/m} .$$
In other words, $k D_0 / m$ is stochastically dominated by an
exponential of mean 1, leading to $\E e^{\lambda k D_0 / m}
\leq 1 / (1 - \lambda)$.  Thus by~(\ref{eq:NA}), $\E \exp (\lambda
k W / m) \leq (1 - \lambda)^{-|B|}$, proving the lemma.   \qed
\endproof
\begin{lem} \label{lem:LD}
Let $|B| = j$ and let $W = \sum_{j \in B} D_j$ as in the
previous lemma.  Then
\begin{equation} \label{eq:LD}
\P ({W \over m} \geq {t \over k}) \leq e^{-t} \left ( {e t \over
   j} \right )^j .
\end{equation}
\end{lem}

\noindent
\proof{Proof} Use Markov's inequality
$$\P ({W \over m} \geq {t \over k}) \leq {\E e^{\lambda k W / m}
   \over e^{\lambda t}} .$$
Set $\lambda = 1 - j/t$ and use the previous lemma to get
\begin{eqnarray*}
\P ({W \over m} \geq {t \over k}) & \leq & \left ( 1 - \lambda
\right)^{-j} e^{-\lambda t} \\[2ex]
& = & \left ( {t \over j} \right)^j e^{j - t} ,
\end{eqnarray*}
proving the lemma.    \qed
\endproof

\noindent
\proof{Proof of Lemma~\ref{lem:1}} We can assume without loss of
generality that $j := \delta k$ is an integer and that $n := \ep
m$ is an integer.  By exchangeability, $p(k , m , \ep , \delta)$
is at most ${k+1 \choose j}$ times the probability that $W(B) / m
\geq \ep$ for any particular $B$ of cardinality $j$.  Setting $t =
k \ep$ and plugging in the result of Lemma~\ref{lem:LD} then gives
\begin{eqnarray*}
p(k , m , \ep , \delta) & \leq & {k+1 \choose j} \left ( {\ep k
\over j} \right )^j e^{j - \ep k} \\[2ex]
   & = & {k \choose \delta k} \left ( {\ep \over \delta}\right )^{\delta k}
   e^{(\delta - \ep)k} .
\end{eqnarray*}
The inequality ${a \choose b} \leq (a/b)^b (a/(a-b))^{a-b}$ holds
for all integers $a \geq b \geq 0$ (with $0^0 := 1$) and leads to
the right-hand side of the previous equation being bounded above by
$$\left ( {1 \over \delta} \right )^{\delta k}
   \left ( {1 \over 1 - \delta} \right )^{(1 - \delta) k}
   \left ( {\ep \over \delta} \right)^{\delta k}
   e^{(\delta - \ep)k} .$$
Hence $p(k , m , \ep , \delta) \leq e^{k r(\ep , \delta)}$ where
$$r (\ep , \delta) = \delta (\log \ep - 2 \log \delta +
   \log (1 - \delta)) - \log (1 - \delta) + \delta - \ep .$$
Since $\log \ep$ and $\log (1 - \delta)$ are negative, we have
$$r (\ep , \delta) \leq 2 \delta \log (1/\delta) - \ep + \delta +
   \log (1 /(1 - \delta)) .$$
For sufficiently small $\ep$, hence small $\delta$, we have
$\delta + \log (1 / (1 - \delta)) < (1/2) \delta \log (1/\delta)$,
hence 
$$
r(\ep , \delta) < (5/2) \delta \log (1/\delta) - \ep \leq \ep/2 - 
\ep = -{\ep \over 2} \,,
$$
by the choice of $\delta$.  This finishes the proof.   \qed
\endproof
\section{Concluding Remarks.} \label{sec:conclusion}

\noindent{$\bullet$} A Markov chain $\Gamma$ with state-space ${\cal X}$ and transition
kernel $P$ is {\em transitive} if, for each pair of states $x,y \in
{\cal X}$, there is an invertible mapping $\Phi:{\cal X} \rightarrow
{\cal X}$ so that $\Phi(x)=y$, and
$P(y,\Phi(z)) = P(x,z)$ for all $z \in {\cal X}$.  Random walks, for
example, are transitive Markov chains.  When the underlying Markov
chain $\Gamma$ is transitive, our model has an equivalent percolation
description.  Indeed, given the sample path $\{\Gamma_n\}$, connect two
vertices $m,\ell \in \Z^{+}$ iff
$$\Gamma_m = \Gamma_\ell\, , \;\mbox{but } \; \Gamma_j \neq \Gamma_m
\, \mbox{for } \, m < j < \ell.$$
A coin is chosen for each cluster
(connected component), and labels are generated at each $x \in \Z^{+}$
by flipping this coin.  The coin used for vertices in the cluster of the
origin is
$\bias$-biased, while the coin used in all other
clusters is fair.  The bonds are hidden from an observer, who
must decide which coin was used for the cluster of the origin.  For
certain $\Gamma$ (e.g., for the random walks considered in Section
\ref{sec:counterexample}), there is a phase transition: for $\theta$
sufficiently small, it cannot be determined which coin was used for the
cluster of the origin, while for $\theta$ large enough, the viewer can
distinguish.  This is an example of a $1$-dimensional, long-range,
dependent percolation model which exhibits a phase transition.
Other $1$-dimensional models that exhibit a phase transition
were studied by Aizenman, Chayes, Chayes, and Newman in \cite{ACCN}.

\noindent{$\bullet$} 
In Sections \ref{sec:counterexample} and \ref{sec:example}, we
constructed explicitly renewal processes whose renewal probabilities 
$\{u_n\}$ have prescribed asymptotics.  Alternatively, we could invoke
the following general result.

\noindent
{\sc Kaluza's Theorem }~\cite{Kal}.  \em If $u(0)=1$ and $u(k-1)u(k+1) \geq
u^2(k)$ for $k \geq 1$, then $\{u_k\}$ is a renewal sequence. \rm

\noindent
See \cite{Kal} or \cite[Theorem 5.3.2]{A} for a proof,
and \cite{Lig} for a generalization.

\noindent{$\bullet$} An extended version of the random coin tossing model, 
when the underlying Markov chain is simple random walk on
$\Z$, is studied in \cite{Lev}.  Each vertex $z \in \Z$ is assigned
a coin with bias $\theta(z)$.  At each move of a random walk on $\Z$,
the coin attached to the walk's position is tossed.  In \cite{Lev}, it
is shown that if $|\{z: \theta(z) \neq 0\}|$ is finite, then the
biases $\theta(z)$ can be recovered up to a symmetry of $\Z$.

\noindent {\bf Some unsolved problems}.  
Recall that $\renp$ and
$\renp'$ denote two independent and identically distributed
renewal processes, and $u_n =
\P[\renp_n = 1]$. The distribution of the sequence of coin
tosses, when a coin with bias $\bias$ is used at renewal times, is
denoted by $\mu_\bias$.
\begin{enumerate}
%

 \item Is the quenched moment generating function criterion in Lemma
 \ref{lem:qmgf} sharp?  That is, does $\E [r^
 {\sum_{n=0}^\infty \renp_n \renp'_n}
 \given \renp ] = \infty$ for some $r < 1+\theta^2$
 imply that $\mu_\bias \perp \mu_0$?
\item Does $\mu_{\theta_1} \perp \mu_0$ imply that $\mu_{\theta_1} \perp
\mu_{\theta_2}$ for all $\theta_2 \neq \theta_1$?
\item For renewal sequences exhibiting a phase transition
  at a critical parameter $\theta_c$,  is $\mu_{\theta_c} \perp \mu_0$?
\end{enumerate}

\section*{Acknowledgments.}
 
We thank A. Dembo, J. Steif, and O. Zeitouni
for useful discussions, and W. Woess and T. Liggett for references.


\Line{\AOPaddress{Dept. of Mathematics\\
196 Auditorium Road, U-9\\
Univ. of Connecticut\\
Storrs, CT 06269-3009\\
{\tt levin@math.uconn.edu}}\hfill
\AOPaddress{Dept. of Mathematics\\
480 Lincoln Drive \\
Univ. of Wisconsin\\
Madison, WI 53706\\
{\tt pemantle@math.wisc.edu}} 
\hfill
\AOPaddress{
Mathematics Institute \\
The Hebrew University \\
Givat-Ram \\
Jerusalem 91904, Israel\\
{\tt peres@math.huji.ac.il} \\
and \\
Dept. of Statistics \\
367 Evans Hall \#3860 \\
University of California \\
Berkeley, CA 94720-3860}
}

\end{document}